\documentclass{amsart}

\usepackage{amsmath,amsthm}
\usepackage{amsfonts}
\usepackage{amssymb}
\usepackage{graphicx}
\usepackage{cite}
\usepackage{times}

\theoremstyle{plain}
\newtheorem{theorem}{Theorem}
\newtheorem{lemma}[theorem]{Lemma}
\newtheorem{proposition}[theorem]{Proposition}

\theoremstyle{definition}
\newtheorem{definition}[theorem]{Definition}

\theoremstyle{remark}
\newtheorem{remark}[theorem]{Remark}

\title{Classification of $U_q(\mathfrak{sl}_2)$-Module Algebra
Structures on~the Quantum Plane}

\author{Steven Duplij}

\address{Theory Group,
Nuclear Physics Laboratory, V. N. Karazin Kharkov National
University, Svoboda Sq. 4, Kharkov 61077, Ukraine}

\urladdr{http://webusers.physics.umn.edu/\~{}duplij} 

\email{steven.a.duplij@univer.kharkov.ua, sduplij@gmail.com}

\author{Sergey Sinel'shchikov}

\address{Mathematics Division,
B. I. Verkin Institute for Low Temperature Physics and Engineering,
Lenin Ave. 47, National Academy of Sciences of Ukraine, Kharkov
61103, Ukraine} 

\email{sinelshchikov@ilt.kharkov.ua}

\date{May 12, 2009 (Revised August 31, 2010)}

\begin{document}

\protect\maketitle 

\begin{abstract}

A complete list of $U_q(\mathfrak{sl}_2)$-module algebra structures on the
quantum plane is produced and the (uncountable family of) isomorphism
classes of these structures are described. The composition series of
representations in question are computed. The classical limits of the
$U_q(\mathfrak{sl}_2)$-module algebra structures are discussed. \vskip2mm

\smallskip


\end{abstract}\smallskip


\begin{center}
\section{Introduction}
\end{center}

The quantum plane \cite{manin2} is known to be a starting point in studying
modules over quantum universal enveloping algebras \cite{Drinf2_engl}. The
structures existing on the quantum plane are widely used as a background to
produce associated structures for more sophisticated quantum algebras
\cite{dup/sin1,dup/li3,dup/li2}. There is one distinguished structure of
$U_q(\mathfrak{sl}_2)$-module algebra on the quantum plane which was widely
considered before (see, e.g., \cite{kassel}). In addition, one could
certainly mention the structure $h(\mathbf{v})=\varepsilon(h)\mathbf{v}$,
where $h\in U_q(\mathfrak{sl}_2)$, $\varepsilon$ is the counit,
$\mathbf{v}$ is a polynomial on the  quantum \ plane. \ Normally \ it \ is
\ disregarded \ because \ of its \ triviality. \ Nevertheless, it turns out
that there exist more (in fact, an uncountable family of~nonisomorphic)
$U_q(\mathfrak{sl}_2)$-module algebra structures which are nontrivial and
can be used in further development of the quantum group theory.

In this paper we suggest a complete description and classification of
$U_q(\mathfrak{sl}_2)$-module algebra structures  existing on the quantum
plane. Specifically, in Section~3 
we use a general form of the automorphism of quantum plane to render the
notion of weight for $U_q(\mathfrak{sl}_2)$-actions considered here. In Section~4 
we present our classification in terms of a pair of symbolic matrices,
which relies upon considering the low dimensional ($0$-th and $1$-st)
homogeneous components of an action. In Section~5 
we describe the composition series for the above structures viewed as
representations in vector spaces.
\bigskip

{\centering\section{Preliminaries}}
\smallskip

Let $H$ be a Hopf algebra whose comultiplication is $\Delta$, counit is
$\varepsilon$, and antipode is $S$ \cite{abe}. Also let $A$ be a unital
algebra with unit $\mathbf{1}$. We will also use the Sweedler notation
$\Delta\left( h\right) =\sum_{i}h_{i}^{\prime}\otimes h_{i}^{\prime\prime}$
\cite{sweedler}.\medskip

\begin{definition}
By a structure of $H$-module algebra on $A$ we mean a homomorphism
$\pi:H\to\operatorname{End}_\mathbb{C}A$ such that:

(i) $\pi(h)(ab)=\sum_i\pi(h_i')(a)\cdot\pi(h_i'')(b)$ for all $h\in H$,
$a,b\in A$;

(ii) $\pi(h)(\mathbf{1})=\varepsilon(h)\mathbf{1}$ for all $h\in H$.

The structures $\pi_1,\pi_2$ are said to be isomorphic if there exists an
automorphism $\Psi$ of the algebra $A$ such that
$\Psi\pi_1(h)\Psi^{-1}=\pi_2(h)$ for all $h\in H$.
\end{definition}

Throughout the paper we assume that $q\in\mathbb{C}\setminus\{0\}$
is not a root of the unit ($q^n\ne 1$ for all non-zero integers
$n$). Consider the quantum plane which is a unital algebra
$\mathbb{C}_q[x,y]$ with two generators $x,y$ and a single
relation
\begin{equation}\label{1}
yx=qxy.
\end{equation}

The quantum universal enveloping algebra $U_{q}\left(  \mathfrak{sl}%
_{2}\right) $ is a unital associative algebra determined by its (Chevalley)
generators $\mathsf{k}$, $\mathsf{k}^{-1}$, $\mathsf{e}$, $\mathsf{f}$, and
the relations%
\begin{align}
\mathsf{k}^{-1}\mathsf{k}  &  =\mathbf{1},\ \ \ \ \mathsf{kk}^{-1}%
=\mathbf{1},\label{kk1}\\
\mathsf{ke}  &  =q^{2}\mathsf{ek},\label{ke}\\
\mathsf{kf}  &  =q^{-2}\mathsf{fk},\label{kf}\\
\mathsf{ef}-\mathsf{fe}  &  =\dfrac{\mathsf{k}-\mathsf{k}^{-1}}{q-q^{-1}}.
\label{effe}%
\end{align}

The standard Hopf algebra structure on $U_q(sl_2)$ is determined by
\begin{align}
\Delta(\mathsf{k}) &=\mathsf{k}\otimes\mathsf{k}, & \label{k0}\\
\Delta(\mathsf{e}) &=\mathbf{1}\otimes\mathsf{e}+
\mathsf{e}\otimes\mathsf{k}, & \label{definition}\\
\Delta(\mathsf{f}) &=\mathsf{f}\otimes\mathbf{1}+
\mathsf{k}^{-1}\otimes\mathsf{f}, & \label{def1}\\
\mathsf{S}(\mathsf{k}) &=\mathsf{k}^{-1}, &
\mathsf{S}(\mathsf{e}) &=-\mathsf{ek}^{-1}, &
\mathsf{S}(\mathsf{f}) &=-\mathsf{kf}, & \nonumber\\
\varepsilon(\mathsf{k}) &=1, &
\varepsilon(\mathsf{e}) &=\varepsilon(\mathsf{f})=0. &&\nonumber
\end{align}

{\center\section{\label{Auto}Automorphisms of the Quantum Plane}}
\smallskip

Denote by $\mathbb{C}_q[x,y]_i$ the $i$-th homogeneous component of
$\mathbb{C}_q[x,y]$, which is a linear span of the monomials $x^my^n$ with
$m+n=i$. Also, given a polynomial $p\in\mathbb{C}_q[x,y]$, denote by
$(p)_i$ the $i$-th homogeneous component of $p$, that is the projection of
$p$ onto $\mathbb{C}_q[x,y]_i$ parallel to the direct sum of all other
homogeneous components of $\mathbb{C}_q[x,y]$.

We rely upon a result by J. Alev and M. Chamarie which gives, in
particular, a description of automorphisms of the algebra
$\mathbb{C}_q[x,y]$ \cite[Prop.~1.4.4(i)]{ale/cha}. In fact, their claim is
much more general, so in the special case we need here we present a quite
elementary proof for the reader's convenience.\medskip

\begin{proposition}\label{T1}
Let $\Psi$ be an automorphism of $\mathbb{C}_q[x,y]$, then there exist
nonzero constants $\alpha,\beta$ such that
\begin{equation}\label{psi}
\Psi:x\mapsto\alpha x,\qquad y\mapsto\beta y.
\end{equation}
\end{proposition}\smallskip

First note that an automorphism as in (\ref{psi}) is well defined on the
entire algebra, because the ideal of relations generated by (\ref{1}) is
$\Psi$-invariant. We split the proof into a series of
lemmas.

\begin{lemma}\label{L1}
One has $(\Psi(x))_0=(\Psi(y))_0=0$.
\end{lemma}

P r o o f. \ We start with proving $(\Psi(x))_0=0$. Suppose the contrary,
that is $(\Psi(x))_0\ne 0$. As $\Psi(y)\ne 0$, we choose the lowest $i$
with $(\Psi(y))_i\ne 0$. Apply $\Psi$ to the relation $yx=qxy$ and then
project it to the $i$-th homogeneous component of $\mathbb{C}_q[x,y]$
(parallel to the direct sum of all other homogeneous components) to get
$(\Psi(y)\Psi(x))_i=q(\Psi(x)\Psi(y))_i$. Clearly, $(\Psi(y)\Psi(x))_i$ is
the lowest homogeneous component of $\Psi(y)\Psi(x)$, and
$(\Psi(y)\Psi(x))_i=(\Psi(y))_i(\Psi(x))_0$. In a similar way
$q(\Psi(x)\Psi(y))_i=q(\Psi(x))_0(\Psi(y))_i$. Because $(\Psi(x))_0$ is a
constant, it commutes with $(\Psi(y))_i$, then
$(\Psi(y))_i(\Psi(x))_0=q(\Psi(y))_i(\Psi(x))_0$, and since $(\Psi(x))_0\ne
0$, we also have $(\Psi(y))_i=q(\Psi(y))_i$. Recall that $q\ne 1$, hence
$(\Psi(y))_i=0$ which contradicts to our choice of $i$. Thus our claim is
proved. The proof of another claim goes in a similar way.
\hfill\rule{0.5em}{0.5em}\smallskip

\begin{lemma}\label{L2}
One has $(\Psi(x))_1\ne 0$, $(\Psi(y))_1\ne 0$.
\end{lemma}\smallskip

P r o o f. \ Let us prove that $(\Psi(x))_1\ne 0$. Suppose the contrary,
which by virtue of Lemma \ref{L1} means that
$\Psi(x)=\sum_ia_ix^{m_i}y^{n_i}$ with $m_i+n_i>1$. The subsequent
application of the inverse automorphism gives $\Psi^{-1}(\Psi(x))$ which is
certainly $x$. On the other hand,
$$\Psi^{-1}(\Psi(x))=\sum_ia_i(\Psi^{-1}(x))^{m_i}(\Psi^{-1}(y))^{n_i}.$$
By Lemma \ref{L1} every nonzero monomial in $\Psi^{-1}(x)$ and
$\Psi^{-1}(y)$ has degree at least one, which implies that
$\Psi^{-1}(\Psi(x))$ is a sum of monomials of degree at least $2$. In
particular, $\Psi^{-1}(\Psi(x))$ can not be $x$. This contradiction proves
the claim. The rest of the statements can be proved in a similar
way.\hfill\rule{0.5em}{0.5em}\smallskip

\begin{lemma}\label{L3}
There exist nonzero constants $\alpha,\beta,\gamma,\delta$ such that
$(\Psi(x))_1=\alpha x$, $(\Psi(y))_1=\beta y$.
\end{lemma}\smallskip

P r o o f. \ Let us apply $\Psi$ to \eqref{1}, then project it to
$\mathbb{C}_q[x,y]_2$ to get $(\Psi(y)\Psi(x))_2=q(\Psi(x)\Psi(y))_2$. It
follows from Lemmas \ref{L1}, \ref{L2} that
$(\Psi(y)\Psi(x))_2=(\Psi(y))_1(\Psi(x))_1$ and
$(\Psi(x)\Psi(y))_2=(\Psi(x))_1(\Psi(y))_1$. Let $(\Psi(x))_1=\alpha x+\mu
y$ and $(\Psi(y))_1=\beta y+\nu x$, which leads to $(\beta y+\nu x)(\alpha
x+\mu y)=q(\alpha x+\mu y)(\beta y+\nu x)$. This, together with \eqref{1}
and Lemma \ref{L2}, implies that $\mu=\nu=0$, $\alpha\ne 0$, and $\beta\ne
0$. \hfill\rule{0.5em}{0.5em}\smallskip

\medskip

Denote by $\mathbb{C}[x]$ and $\mathbb{C}[y]$ the linear spans of
$\{x^n|\:n\ge 0\}$ and $\{y^n|\:n\ge 0\}$, respectively.
Obviously, one has the direct sum decompositions
$$
\mathbb{C}_q[x,y]=\mathbb{C}[x]\oplus y\mathbb{C}_q[x,y]=
\mathbb{C}[y]\oplus x\mathbb{C}_q[x,y].
$$

Given any polynomial $P\in\mathbb{C}_{q}[x,y]$, let $\left( P\right) _{x}$
be its projection to $\mathbb{C}[x]$ parallel to $y\mathbb{C}_{q}[x,y]$,
and in a similar way define $\left( P\right) _{y}$. Obviously,
$\mathbb{C}[x]$ and $\mathbb{C}[y]$ are commutative subalgebras.\medskip

\begin{lemma}\label{L4}
One has $(\Psi(x))_y=(\Psi(y))_x=0$.
\end{lemma}\smallskip

P r o o f. \ First we prove that $(\Psi(x))_y=0$. Project $yx=qxy$ to
$\mathbb{C}[y]$ to obtain $(\Psi(y))_y(\Psi(x))_y=q(\Psi(x))_y(\Psi(y))_y$.
On the other hand, $(\Psi(y))_y(\Psi(x))_y=(\Psi(x))_y(\Psi(y))_y$, so that
$(1-q)(\Psi(x))_y(\Psi(y))_y=0$. Since $q\ne 1$, we deduce that
$(\Psi(x))_y(\Psi(y))_y=0$. It follows from Lemma \ref{L3} that
$(\Psi(y))_y\ne 0$, and since $\mathbb{C}_q[x,y]$ is a domain
\cite{jantzen}, we finally obtain $(\Psi(x))_y=0$. The proof of another
claim goes in a similar way. \hfill\rule{0.5em}{0.5em}

\medskip

P r o o f of Proposition \ref{T1}. \ It follows from Lemma \ref{L4} that
$\Psi(x)=xP$ for some $P\in\mathbb{C}_q[x,y]$. An application of
$\Psi^{-1}$ gives $x=\Psi^{-1}(x)\Psi^{-1}(P)$. Since $\deg x=1$, one
should have either $\deg\Psi^{-1}(x)=0$ or $\deg\Psi^{-1}(P)=0$.
Lemma~\ref{L1} implies that $\deg\Psi^{-1}(x)\ne 0$, hence
$\deg\Psi^{-1}(P)=0$, that is $\Psi^{-1}(P)$ is a nonzero constant, \ and \
so \ $P=\Psi\Psi^{-1}(P)$ \ is \ the \ same \ constant \ (we denote it
by~$\alpha$). The second claim can be proved in a similar
way.\hfill\rule{0.5em}{0.5em}

{\center\section{\label{Str}\boldmath The Structures of
$U_q(\mathfrak{sl}_2)$-Module Algebra  on the Quantum Plane}}

We describe here the $U_{q}\left( \mathfrak{sl}_{2}\right) $-module algebra
structures on $\mathbb{C}_{q}[x,y]$ and then classify them up to
isomorphism.

For the sake of brevity, given a $U_q(\mathfrak{sl}_2)$-module algebra
structure on $\mathbb{C}_q[x,y]$, we can associate a $2\times 3$ matrix
with entries from $\mathbb{C}_q[x,y]$
\begin{equation}
\mathsf{M}\overset{definition}{=}\left\Vert
\begin{array}[c]{c}
\mathsf{k}\\
\mathsf{e}\\
\mathsf{f}
\end{array}
\right\Vert\cdot\left\Vert x,y\right\Vert=\left\Vert
\begin{array}[c]{cc}
\mathsf{k}(x) & \mathsf{k}(y) \\
\mathsf{e}(x) & \mathsf{e}(y) \\
\mathsf{f}(x) & \mathsf{f}(y)
\end{array}
\right\Vert , \label{m}
\end{equation}
where $\mathsf{k}$, $\mathsf{e}$, $\mathsf{f}$ are the generators of
$U_q(\mathfrak{sl}_2)$ and $x$, $y$ are the generators of
$\mathbb{C}_q[x,y]$. We call $\mathsf{M}$ a \textit{full action matrix}.
Conversely, suppose we have a matrix $\mathsf{M}$ with entries from
$\mathbb{C}_q[x,y]$ as in \eqref{m}. To derive the associated
$U_q(\mathfrak{sl}_2)$-module algebra structure on $\mathbb{C}_q[x,y]$ we
set (using the Sweedler notation)
\begin{align}
(\mathsf{ab})u\overset{definition}{=}\mathsf{a}(\mathsf{b}u),\qquad
\mathsf{a},\mathsf{b} & \in U_q(\mathfrak{sl}_2),\qquad
u\in\mathbb{C}_q[x,y],\label{abu}\\
\mathsf{a}(uv)\overset{definition}{=}
\Sigma_i(\mathsf{a}_i'u)\cdot(\mathsf{a}_i''v),\qquad\mathsf{a} & \in
U_q(\mathfrak{sl}_2),\qquad u,v\in\mathbb{C}_q[x,y], \label{auv}
\end{align}
which determines a well-defined action of $U_q(\mathfrak{sl}_2)$ on
$\mathbb{C}_q[x,y]$ iff the following properties hold. Firstly, an
application (defined by \eqref{abu}) of an element from the relation ideal
of $U_q(\mathfrak{sl}_2)$ \eqref{kk1}--\eqref{effe} to any
$u\in\mathbb{C}_q[x,y]$ should produce zero. Secondly, a result of
application (defined by \eqref{auv}) of any $\mathsf{a}\in
U_q(\mathfrak{sl}_2)$ to an element of the relation ideal of
$\mathbb{C}_q[x,y]$ \eqref{1} vanishes. These conditions are to be verified
in the specific cases considered below.

Note that, given a $U_{q}\left( \mathfrak{sl}_{2}\right) $-module
algebra structure on the quantum plane, the action of the
generator $\mathsf{k}$ determines an automorphism of
$\mathbb{C}_{q}[x,y]$, which is a consequence of invertibility of
$\mathsf{k}$ and $\Delta\left( \mathsf{k}\right)
=\mathsf{k}\otimes\mathsf{k}$. In particular, it follows from
(\ref{psi}) that $\mathsf{k}$ is determined completely by its
action $\Psi$ on the generators
presented by a $1\times2$-matrix $\mathsf{M}_{\mathsf{k}}$ as follows%
\begin{equation}
\mathsf{M}_{\mathsf{k}}\overset{definition}{=}\left\Vert \mathsf{k}\left(  x\right)
,\mathsf{k}\left(  y\right)  \right\Vert =\left\Vert \alpha x,\beta
y\right\Vert \label{mk}%
\end{equation}
for some $\alpha,\beta\in\mathbb{C\setminus}\left\{ 0\right\} $(which is
certainly a minor of $\mathsf{M}$ (\ref{m})). Therefore every monomial
$x^{n}y^{m}\in\mathbb{C}_{q}[x,y]$ is an eigenvector for $\mathsf{k}$, and
the associated eigenvalue $\alpha^{n}\beta^{m}$ will be referred to as a
\textit{weight} of this monomial, which will be written as
$\mathbf{wt}\left( x^{n}y^{m}\right) =\alpha^{n}\beta^{m}$.

We will also need another minor of $\mathsf{M}$ as follows
\begin{equation}
\mathsf{M}_{\mathsf{ef}}\overset{definition}{=}\left\Vert
\begin{array}[c]{cc}
\mathsf{e}(x) & \mathsf{e}(y) \\
\mathsf{f}(x) & \mathsf{f}(y)
\end{array}
\right\Vert, \label{mef}
\end{equation}
and we call $\mathsf{M}_\mathsf{k}$ and $\mathsf{M}_\mathsf{ef}$ an
\textit{action }$\mathsf{k}$-\textit{matrix} and an action $\mathsf{ef}
$-\textit{matrix}, respectively.\newpage

It follows from \eqref{ke}--\eqref{kf} that each entry of $\mathsf{M}$ is a
weight vector, in particular, all the nonzero monomials which constitute a
specific entry should be of the same weight. Specifically, by some abuse of
notation we can write
\begin{align*}
& \mathbf{wt}(\mathsf{M})\overset{definition}{=}\left(
\begin{array}[c]{cc}
\mathbf{wt}(\mathsf{k}(x)) & \mathbf{wt}(\mathsf{k}(y)) \\
\mathbf{wt}(\mathsf{e}(x)) & \mathbf{wt}(\mathsf{e}(y)) \\
\mathbf{wt}(\mathsf{f}(x)) & \mathbf{wt}(\mathbf{\mathsf{f}}(y))
\end{array}\right)\\
& \bowtie\left(
\begin{array}[c]{cc}
\mathbf{wt}(x) & \mathbf{wt}(y) \\
q^2\mathbf{wt}(x) & q^2\mathbf{wt}(y) \\
q^{-2}\mathbf{wt}(x) & q^{-2}\mathbf{wt}(y)
\end{array}
\right)=\left(
\begin{array}[c]{cc}
\alpha & \beta\\
q^2\mathbf{\alpha} & q^2\mathbf{\beta}\\
q^{-2}\mathbf{\alpha} & q^{-2}\mathbf{\beta}
\end{array}
\right),
\end{align*}
where the relation $\bowtie$ between the two matrices $A=(a_{ij})$ and
$B=(b_{ij})$ is defined as follows:\smallskip

{\bf Notation.} $A\bowtie B$ if for every pair of indices $i,j$ such that
both $a_{ij}$ and $b_{ij}$ are nonzero, one has $a_{ij}=b_{ij}$, e.g.,
$\begin{pmatrix}1 & 0\\ 0 & 2\end{pmatrix}\bowtie\begin{pmatrix}1 & 3\\ 0 &
0\end{pmatrix}$.\smallskip

As an immediate consequence, we also have\smallskip

\begin{proposition}\label{Pmono}
Suppose that $\alpha$/$\beta$ is not a root of the unit. Then every
homogeneous component $(\mathsf{e}(x))_n$, $(\mathsf{e}(y))_n$,
$(\mathsf{f}(x))_n$, $(\mathsf{f}(y))_n$, $n\ge 0$, if nonzero, reduces to
a monomial.
\end{proposition}

P r o o f. Under our assumptions on $\alpha$, $\beta$, the weights of the
monomials $x^{i}y^{n-i}$, $0\leq i\leq n$, of degree $n$ are pairwise
different. Since $\mathsf{e}(x)$, $\mathsf{e}(y)$, $\mathsf{f}(x)$,
$\mathsf{f}(y)$ are weight vectors, our claim follows.
\hfill\rule{0.5em}{0.5em}\smallskip 

Our basic observation is that the $U_q(\mathfrak{sl}_2)$-actions in
question are actually determined to a large extent by the projections of
$\mathsf{M}$ to the lower homogeneous components of $\mathbb{C}_q[x,y]$.

Next, we denote by $(M)_i$ the $i$-th homogeneous component of $M$, whose
elements are just the $i$-th homogeneous components of the corresponding
entries of $M$. Thus every matrix element of $M$, if nonzero, admits a
well-defined weight.

Let us introduce the constants $a_0,b_0,c_0,d_0\in\mathbb{C}$ such that
zero degree component of the full action matrix is
\begin{equation}
(\mathsf{M})_0=\left(
\begin{array}[c]{cc}
0 & 0\\
a_0 & b_0\\
c_0 & d_0
\end{array}
\right)_0. \label{m0}
\end{equation}
Here we keep the subscript $0$ to the matrix in the r.h.s. to emphasize the
origin of this matrix as the $0$-th homogeneous component of $\mathsf{M}$.
Note that the weights of nonzero projections of (weight) entries of
$\mathsf{M}$ should have the same weight. Hence
\begin{equation}
\mathbf{wt}\left((\mathsf{M})_0\right)\bowtie\left(
\begin{array}[c]{cc}
0 & 0\\
q^2\mathbf{\alpha} & q^2\mathbf{\beta}\\
q^{-2}\mathbf{\alpha} & q^{-2}\mathbf{\beta}
\end{array}
\right)_0. \label{wtm0}
\end{equation}
On the other hand, as all the entries of $(\mathsf{M})_0$ are constants
\eqref{m0}, one also deduces
\begin{equation}
\mathbf{wt}\left((\mathsf{M})_0\right)\bowtie\left(
\begin{array}[c]{cc}
0 & 0\\
1 & 1\\
1 & 1
\end{array}
\right)_0, \label{wtm11}
\end{equation}
where the relation $\bowtie$ is understood as a set of elementwise
equalities iff they are applicable, that is, when the corresponding entry
of the projected matrix $(\mathsf{M})_0$ is nonzero. Therefore, it is not
possible to have all nonzero entries in the $0$-th homogeneous component of
$\mathsf{M}$ simultaneously.

The classification of $U_q(\mathfrak{sl}_2)$-module algebra structures on
the quantum plane we are about to suggest will be done in terms of a pair
of symbolic matrices derived from the minor $\mathsf{M}_\mathsf{ef}$ only.
Now we use $(\mathsf{M}_\mathsf{ef})_i$ to construct a symbolic matrix
$\left(\overset{\star}{\mathsf{M}}_\mathsf{ef}\right)_i$ whose entries are
symbols $\mathbf{0}$ or $\mathbf{\star}$ as follows: a nonzero entry of
$(\mathsf{M}_\mathsf{ef})_i$ is replaced by $\star$, while a zero entry is
replaced by the symbol $\mathbf{0}$.

In the case of $0$-th components the specific elementwise
relations involved in (\ref{wtm0}) imply that each column of
$\left(\overset{\star}{\mathsf{M}}_{\mathsf{ef}}\right)_0$ should
contain at least one $\mathbf{0}$, and so
$\left(\overset{\star}{\mathsf{M}}_{\mathsf{ef}}\right)_0$ can be
either of the following 9 matrices
\begin{align}
&  \left(
\begin{array}[c]{cc}
\mathbf{0} & \mathbf{0}\\
\mathbf{0} & \mathbf{0}
\end{array}
\right)_0, \nonumber\\
&  \left(
\begin{array}[c]{cc}
\mathbf{\star} & \mathbf{0}\\
\mathbf{0} & \mathbf{0}
\end{array}
\right)_0,\left(
\begin{array}[c]{cc}
\mathbf{0} & \mathbf{\star}\\
\mathbf{0} & \mathbf{0}
\end{array}
\right)_0,\left(
\begin{array}[c]{cc}
\mathbf{0} & \mathbf{0}\\
\mathbf{\star} & \mathbf{0}
\end{array}
\right)_0,\left(
\begin{array}[c]{cc}
\mathbf{0} & \mathbf{0}\\
\mathbf{0} & \mathbf{\star}
\end{array}
\right)_0, \label{0}\\
&  \left(
\begin{array}[c]{cc}
\mathbf{\star} & \mathbf{\star}\\
\mathbf{0} & \mathbf{0}
\end{array}
\right)_0,\left(
\begin{array}[c]{cc}
\mathbf{0} & \mathbf{0}\\
\mathbf{\star} & \mathbf{\star}
\end{array}
\right)_0,\left(
\begin{array}[c]{cc}
\mathbf{\star} & \mathbf{0}\\
\mathbf{0} & \mathbf{\star}
\end{array}
\right)_0,\left(
\begin{array}[c]{cc}
\mathbf{0} & \mathbf{\star}\\
\mathbf{\star} & \mathbf{0}
\end{array}
\right)_0. \nonumber
\end{align}

An application of $\mathsf{e}$ and $\mathsf{f}$ to \eqref{1} by using
\eqref{mk} gives
\begin{align}
y\mathsf{e}(x)-q\beta\mathsf{e}(x)y &=qx\mathsf{e}(y)-\alpha\mathsf{e}(y)x,
\label{exy}\\
\mathsf{f}(x)y-q^{-1}\beta^{-1}y\mathsf{f}(x)
&=q^{-1}\mathsf{f}(y)x-\alpha^{-1}x\mathsf{f}(y). \label{fxy}
\end{align}

After projecting \eqref{exy}--\eqref{fxy} to $\mathbb{C}_q[x,y]_1$ we
obtain
\begin{align*}
a_0(1-q\beta)y &=b_0(q-\alpha)x,\\
d_0\left(1-q\alpha^{-1}\right)x &=c_0\left(q-\beta^{-1}\right)y,
\end{align*}
which certainly implies
$$
a_0(1-q\beta)=b_0(q-\alpha)=d_0\left(1-q\alpha^{-1}\right)=
c_0\left(q-\beta^{-1}\right)=0.
$$
This determines the weight constants $\alpha$ and $\beta$ as follows:
\begin{align}
a_0 & \ne 0\Longrightarrow\beta=q^{-1}, & \label{b1}\\
b_0 & \ne 0\Longrightarrow\alpha=q, & \label{b2}\\
c_0 & \ne 0\Longrightarrow\beta=q^{-1}, & \label{b3}\\
d_0 & \ne 0\Longrightarrow\alpha=q. & \label{b4}
\end{align}

The deduction compared to \eqref{wtm0}, \eqref{wtm11} implies that the
symbolic matrices from \eqref{0} containing two $\star$'s should be
excluded. Also, using \eqref{wtm0} and \eqref{b1}--\eqref{b4} we conclude
that the position of $\star$ in the remaining symbolic matrices completely
determines the associated weight constants by
\begin{align}
\left(
\begin{array}[c]{cc}
\mathbf{\star} & \mathbf{0}\\
\mathbf{0} & \mathbf{0}
\end{array}
\right)_0 && \Longrightarrow && \alpha &=q^{-2}, & \beta &=q^{-1},\label{01}\\
\left(
\begin{array}[c]{cc}
\mathbf{0} & \mathbf{\star}\\
\mathbf{0} & \mathbf{0}
\end{array}
\right)_0 && \Longrightarrow && \alpha &=q, & \beta &=q^{-2},\label{02}\\
\left(
\begin{array}[c]{cc}
\mathbf{0} & \mathbf{0}\\
\mathbf{\star} & \mathbf{0}
\end{array}
\right)_0 && \Longrightarrow && \alpha &=q^2, & \beta &=q^{-1},\label{03}\\
\left(
\begin{array}[c]{cc}
\mathbf{0} & \mathbf{0}\\
\mathbf{0} & \mathbf{\star}
\end{array}
\right)_0 && \Longrightarrow && \alpha &=q, & \beta &=q^2. \label{04}
\end{align}
As for the matrix $\left(
\begin{array}
[c]{cc}%
\mathbf{0} & \mathbf{0}\\
\mathbf{0} & \mathbf{0}%
\end{array}
\right) _{0}$, it does not determine the weight constants at all.

Next, for the $1$-st homogeneous component, one has
$\mathbf{wt}(\mathsf{e}(x))=q^2\mathbf{wt}(x)\ne\mathbf{wt}(x)$ (because
$q^2\ne 1$), which implies $(\mathsf{e}(x))_1=a_1y$, and in a similar way
we have
$$
(\mathsf{M}_{\mathsf{ef}})_1=\left(
\begin{array}[c]{cc}
a_1y & b_1x\\
c_1y & d_1x
\end{array}
\right)_1
$$
with $a_1,b_1,c_1,d_1\in\mathbb{C}$. This allows us to introduce a symbolic
matrix $\left(\overset{\star}{\mathsf{M}}_{\mathsf{ef}}\right)_1$ as above.
Using the relations between the weights similar to \eqref{wtm0}, we obtain
\begin{equation}
\mathbf{wt}((\mathsf{M}_\mathsf{ef})_1)\bowtie\left(
\begin{array}[c]{cc}
q^2\alpha & q^2\beta\\
q^{-2}\alpha & q^{-2}\beta
\end{array}
\right)_1\bowtie\left(
\begin{array}[c]{cc}
\beta & \alpha\\
\beta & \alpha
\end{array}
\right)_1, \label{wtm1}
\end{equation}
here $\bowtie$ is implicit for a set of the elementwise equalities
applicable iff the respective entry of the projected matrix
$(\mathsf{M})_1$ is nonvanishing.

This means that every row and every column of $\left( \overset{\star
}{\mathsf{M}}_{\mathsf{ef}}\right)_1$ may contain at least one
$\mathbf{0}$. Now project \eqref{exy}--\eqref{fxy} to $\mathbb{C}_q[x,y]_2$
to obtain
\begin{align*}
a_1(1-q\beta)y^2 &=b_1(q-\alpha)x^2,\\
d_1\left(1-q\alpha^{-1}\right)x^2 &=c_1\left(q-\beta^{-1}\right)y^2,
\end{align*}
whence $a_1(1-q\beta)=b_1(q-\alpha)=d_1\left(1-q\alpha^{-1}\right)
=c_1\left(q-\beta^{-1}\right)=0$. As a consequence we have
\begin{align}
a_1 & \ne 0 & \Longrightarrow && \beta &=q^{-1}, & \label{aa1}\\
b_1 & \ne 0 & \Longrightarrow && \alpha &=q, & \label{aa2}\\
c_1 & \ne 0 & \Longrightarrow && \beta &=q^{-1}, & \label{aa3}\\
d_1 & \ne 0 & \Longrightarrow && \alpha &=q. & \label{aa4}
\end{align}

A comparison of \eqref{wtm1} with \eqref{aa1}--\eqref{aa4} allows one to
discard the symbolic matrix $\left(
\begin{array}[c]{cc}
\mathbf{\star} & \mathbf{0}\\
\mathbf{0} & \mathbf{\star}
\end{array}
\right)_1$ from the list of symbolic matrices with at least one
$\mathbf{0}$ at every row or column. As for other symbolic
matrices with the above property, we get
\begin{align}
\left(
\begin{array}[c]{cc}
\mathbf{\star} & \mathbf{0}\\
\mathbf{0} & \mathbf{0}
\end{array}
\right)_1 & \Longrightarrow\alpha=q^{-3},\quad\beta=q^{-1},\label{11}\\
\left(
\begin{array}[c]{cc}
\mathbf{0} & \mathbf{\star}\\
\mathbf{0} & \mathbf{0}
\end{array}
\right)_1 & \Longrightarrow\alpha=q,\quad\beta=q^{-1},\label{12}\\
\left(
\begin{array}[c]{cc}
\mathbf{0} & \mathbf{0}\\
\mathbf{\star} & \mathbf{0}
\end{array}
\right)_1 & \Longrightarrow\alpha=q,\quad\beta=q^{-1},\label{13}\\
\left(
\begin{array}[c]{cc}
\mathbf{0} & \mathbf{0}\\
\mathbf{0} & \mathbf{\star}
\end{array}
\right)_1 & \Longrightarrow\alpha=q,\quad\beta=q^3,\label{14}\\
\left(
\begin{array}[c]{cc}
\mathbf{0} & \mathbf{\star}\\
\mathbf{\star} & \mathbf{0}
\end{array}
\right)_1 & \Longrightarrow\alpha=q,\quad\beta=q^{-1}. \label{15}
\end{align}

The matrix $\left(
\begin{array}
[c]{cc}%
\mathbf{0} & \mathbf{0}\\
\mathbf{0} & \mathbf{0}%
\end{array}
\right)_{1}$ does not determine the weight constants in the way described
above.

In view of the above observations we see that in most cases a pair of
symbolic matrices corresponding to $0$-th and $1$-st homogeneous components
determines completely the weight constants of the conjectured associated
actions. It will be clear from the subsequent arguments that the higher
homogeneous components are redundant within the presented classification.
Therefore, we introduce the table of families of
$U_q(\mathfrak{sl}_2)$-module algebra structures, each family is labelled
by two symbolic matrices
$\left(\overset{\star}{\mathsf{M}}_{\mathsf{ef}}\right)_0$,
$\left(\overset{\star}{\mathsf{M}}_{\mathsf{ef}}\right)_1$, and we call
such a family a
$\left[\left(\overset{\star}{\mathsf{M}}_{\mathsf{ef}}\right)_0;
\left(\overset{\star}{\mathsf{M}}_{\mathsf{ef}}\right)_1\right]$-series.
Note that the series labelled with pairs of nonzero symbolic matrices at
both positions are empty, because each of the matrices determines a pair of
specific weight constants $\alpha$ and $\beta$ \eqref{01}--\eqref{04} which
fails to coincide to any pair of such constants associated to the set of
nonzero symbolic matrices at the second position \eqref{11}--\eqref{15}.
Also, the series with zero symbolic matrix at the first position and
symbolic matrices containing only one $\star$ at the second position are
empty.\smallskip

For instance, show that $\left[\left(
\begin{array}{cc}\mathbf{0} & \mathbf{0}\\ \mathbf{0} & \mathbf{0}\end{array}
\right)_0; \left(\begin{array}{cc}\mathbf{\star} & \mathbf{0}\\
\mathbf{0} & \mathbf{0}\end{array}\right)_1\right]$-series is empty. If we
suppose the contrary, then it follows from \eqref{effe} that within this
series we have
$$e(f(x))-f(e(x))=-(1+q^2+q^{-2})x.$$
We claim that the projection of the l.h.s. to $\mathbb{C}_q[x,y]_1$ is
zero. Start with observing that, if the first symbolic matrix consists of
$\mathbf{0}$'s only, one cannot reduce a degree of any monomial by applying
$e$ or $f$. On the other hand, within this series $f(x)$ is a sum of the
monomials whose degree is at least $2$. Therefore, the term $e(f(x))$ has
zero projection to $\mathbb{C}_q[x,y]_1$. Similarly, $f(e(x))$ has also
zero projection to $\mathbb{C}_q[x,y]_1$. The contradiction we get proves
our claim.

In a similar way, one can prove that all other series with zero symbolic
matrix at the first position and symbolic matrices containing only one
$\star$ at the second position are empty.

In the framework of our classification we obtained 24 \textquotedblleft
empty\textquotedblright
$\left[\!\left(\!\overset{\star}{\mathsf{M}}_{\mathsf{ef}}\!\right)_{\!0}\!\!;
\left(\!\overset{\star}{\mathsf{M}}_{\mathsf{ef}}\!\right)_{\!\!1}\right]$-series.
Next turn to \textquotedblleft nonempty\textquotedblright\ series. We start
with the simplest case in which the action $\mathsf{ef}$-matrix is zero,
while the full action matrix is
$$
\mathsf{M}=\left\Vert
\begin{array}
[c]{cc}
\alpha x & \beta y\\
0 & 0\\
0 & 0
\end{array}
\right\Vert .
$$

\begin{theorem}\label{T01}
The $\left[\left(
\begin{array}
[c]{cc}
\mathbf{0} & \mathbf{0}\\
\mathbf{0} & \mathbf{0}
\end{array}
\right)_0;\left(
\begin{array}
[c]{cc}
\mathbf{0} & \mathbf{0}\\
\mathbf{0} & \mathbf{0}
\end{array}
\right)_1\right]$-series consists of 4 $U_q(\mathfrak{sl}_2)$-module
algebra structures on the quantum plane given by
\begin{align}
\mathsf{k}(x) &=\pm x,\qquad\mathsf{k}(y)=\pm y, &\label{t1}\\
\mathsf{e}(x) &=\mathsf{e}(y)=\mathsf{f}(x)=\mathsf{f}(y)=0,& \label{t2}
\end{align}
which are pairwise nonisomorphic.
\end{theorem}

P r o o f. \ It is evident that \eqref{t1}--\eqref{t2} determine a
well-defined $U_q(\mathfrak{sl}_2)$-action consistent with the
multiplication in $U_q(\mathfrak{sl}_2)$ and in the quantum plane, as well
as with comultiplication in $U_q(\mathfrak{sl}_2)$. Prove that there are no
other $U_q(\mathfrak{sl}_2)$-actions here. Note that an application of the
l.h.s. of \eqref{effe} to $x$ or $y$ has zero projection to
$\mathbb{C}_q[x,y]_1$, because in this series $\mathsf{e}$ and $\mathsf{f}$
send any monomial to a sum of the monomials of higher degree. Therefore,
$\left(\mathsf{k}-\mathsf{k}^{-1}\right)(x)=
\left(\mathsf{k}-\mathsf{k}^{-1}\right)(y)=0$, and hence
$\alpha-\alpha^{-1}=\beta-\beta^{-1}=0$, which leads to
$\alpha,\beta\in\{1,-1\}$. To prove \eqref{t2}, note that
$\mathbf{wt}(\mathsf{e}(x))=q^2\mathbf{wt}(x)=\pm q^2\ne\pm 1$. On the
other hand, the weight of any nonzero weight vector in this series is $\pm
1$. This and similar arguments which involve $\mathsf{e}$, $\mathsf{f}$,
$x$, $y$ imply \eqref{t2}.

To see that the $U_q(\mathfrak{sl}_2)$-module algebra structures are
pairwise non-isomorphic, observe that all the automorphisms of the quantum
plane commute with the action of $\mathsf{k}$ (see Sect.~3
).\hfill\rule{0.5em}{0.5em}\smallskip

The action we reproduce in the next theorem is well known \cite{lam/rad,
mon/smi}, and here is the place for it in our
classification.\smallskip

\begin{theorem}
The $\left[\left(
\begin{array}[c]{cc}
\mathbf{0} & \mathbf{0}\\
\mathbf{0} & \mathbf{0}
\end{array}
\right)_0;\left(
\begin{array}[c]{cc}
\mathbf{0} & \mathbf{\star}\\
\mathbf{\star} & \mathbf{0}
\end{array}
\right)_1\right]$-series consists of a one-para\-meter
($\tau\in\mathbb{C}\setminus\{0\}$) family of $U_q(\mathfrak{sl}_2)$-module
algebra structures on the quantum plane
\begin{align}
\mathsf{k}(x) &=qx, & \mathsf{k}(y) &=q^{-1}y,\label{kkqq} &&& \\
\mathsf{e}(x) &=0, & \mathsf{e}(y) &=\tau x,\label{kkqq1} &&& \\
\mathsf{f}(x) &=\tau^{-1}y, & \mathsf{f}(y) &=0. &&& \label{ffxy}
\end{align}

All these structures are isomorphic, in particular, to the action as above
with $\tau=1$.
\end{theorem}

The full action matrix related to \eqref{kkqq}--\eqref{ffxy} is
$$
\mathsf{M}=\left\Vert
\begin{array}[c]{cc}
qx & q^{-1}y\\
0 & x\\
y & 0
\end{array}
\right\Vert .
$$

P r o o f. \ It is easy to check that \eqref{kkqq}--\eqref{ffxy} are
compatible to all the relations in $U_q(\mathfrak{sl}_2)$ and
$\mathbb{C}_q[x,y]$, hence determine a well-defined
$U_q(\mathfrak{sl}_2)$-module algebra structure on the quantum plane
\cite{mon/smi}.\smallskip

Prove that the $\left[\left(
\begin{array}[c]{cc}
\mathbf{0} & \mathbf{0}\\
\mathbf{0} & \mathbf{0}
\end{array}
\right)_0;\left(
\begin{array}[c]{cc}
\mathbf{0} & \mathbf{\star}\\
\mathbf{\star} & \mathbf{0}
\end{array}
\right)_1\right]$-series contains no other actions except those given by
\eqref{kkqq}--\eqref{ffxy}. Let us first prove that the matrix elements of
$\mathsf{M}_{\mathsf{ef}}$ \eqref{mef} contain no terms of degree higher
than one, i.e. $\left(\mathsf{M}_{\mathsf{ef}}\right)_n=0$ for $n\ge 2$. A
general form for $\mathsf{e}(x)$ and $\mathsf{e}(y)$ here is
\begin{equation}\label{er}
\mathsf{e}(x)=\sum_{m+n\ge 2}\bar{\rho}_{mn}x^my^n,\qquad
\mathsf{e}(y)=\tau_\mathsf{e}x+\sum_{m+n\ge 2}\bar{\sigma}_{mn}x^my^n,
\end{equation}
where $\tau_\mathsf{e},\bar{\rho}_{mn},\bar{\sigma}_{mn}\in\mathbb{C}$,
$\tau_\mathsf{e}\ne 0$. Note that in this series
$$
\mathbf{wt}\left(\mathsf{M}_\mathsf{ef}\right)=\left(
\begin{array}[c]{cc}
q^3 & q\\
q^{-1} & q^{-3}
\end{array}
\right) .
$$

In particular, $\mathbf{wt}(\mathsf{e}(x))=q^3$ and
$\mathbf{wt}(\mathsf{e}(y))=q$, which reduces the general form \eqref{er}
to a sum of terms with each one having the same fixed weight
\begin{align}
\mathsf{e}(x) &=\sum_{m\ge 0}\rho_mx^{m+3}y^m,\label{er1}\\
\mathsf{e}(y) &=\tau_\mathsf{e}x+\sum_{m\ge 0}\sigma_mx^{m+2}y^{m+1}.
\label{es1}
\end{align}

Substitute \eqref{er1}--\eqref{es1} to \eqref{exy} and then project it to
the one-dimensional subspace $\mathbb{C}x^{m+3}y^{m+1}$ (for every $m\ge
0$) to obtain
$$\dfrac{\rho_m}{\sigma_m}=-q\dfrac{1-q^{m+1}}{1-q^{m+3}}.$$

In a similar way, the relations
$\mathbf{wt}(\mathsf{f}(x))=q^{-1}$ and
$\mathbf{wt}(\mathsf{f}(y))=q^{-3}$ imply that
\begin{align}
\mathsf{f}(x) &=\tau_fy+\sum_{n\ge 0}\rho_n'x^{n+1}y^{n+2},\label{fr1}\\
\mathsf{f}(y) &=\sum_{n\ge 0}\sigma_n'x^ny^{n+3},\label{fs1}
\end{align}
where $\tau_f\in\mathbb{C}\setminus\{0\}$. An application of
\eqref{fr1}--\eqref{fs1} and \eqref{fxy} with subsequent projection to
$\mathbb{C}x^{n+1}y^{n+3}$ (for every $n\ge 0$) allows one to get
$$\dfrac{\rho_n'}{\sigma_n'}=-q^{-1}\dfrac{1-q^{n+3}}{1-q^{n+1}}.$$

Thus we have
$$
\mathsf{M}_{\mathsf{ef}} =\!\left(\!
\begin{array}[c]{cc}
0 & \tau_ex\\
\tau_fy & 0
\end{array}
\!\right)
 +\sum_{n\ge 0}\left(\!
\begin{array}[c]{cc}
-\mu_nq(1-q^{n+1})x^{n+3}y^n & \mu_n(1-q^{n+3})x^{n+2}y^{n+1}\\
\nu_n(1-q^{n+3})x^{n+1}y^{n+2} & -\nu_nq(1-q^{n+1})x^ny^{n+3}
\end{array}
\!\right)\!,
$$
where $\mu_n,\nu_n\in\mathbb{C}$. We intend to prove that the second matrix
in this sum is zero. Assume the contrary. In the case there exist both
nonzero $\mu_n$'s and $\nu_n$'s, and since the sums here are finite, for
the first row choose the largest index $n_e$ with $\mu_{n_e}\ne 0$ and for
second row, the largest index $n_f$ with $\nu_{n_f}\ne 0$. Then using
\eqref{definition}--\eqref{def1}, we deduce that the highest degree of the
monomials in $(ef-fe)(x)$ is $2n_e+2n_f+5$. This monomial appears to be
unique, and its precise computation gives
$\mu_{n_e}\nu_{n_f}q^{n_en_f-1}(1-q^{n_2+n_f+4})(1-q^{2n_e+2n_f+6})
x^{n_e+n_f+3}y^{n_e+n_f+2}$. Therefore, $(ef-fe)(x)$ has a nonzero
projection onto the one dimensional subspace spanned by the monomial
$x^{n_e+n_f+3}y^{n_e+n_f+2}$, the latter being of degree higher than $1$.
This contradicts to \eqref{effe} whose r.h.s. applied to $x$ has degree
$1$.

In the case when all $\nu_n$'s are zero and some $\mu_n$'s are
nonvanishing we have that the highest degree monomial of
$(ef-fe)(x)$ is of the form
$$
\tau_f\mu_{n_e}\frac{(1-q^{n_e+3})(1-q^{2n_e+4})}{q^{n_e+1}(1-q^2)}
x^{n_e+2}y^{n_e+1},
$$
which is nonzero under our assumptions on $q$. This again produces the same
contradiction as above. In the opposite case when all $\mu_n$'s are zero
and some $\nu_n$'s are nonvanishing, a similar computation works, which
also leads to a contradiction. Therefore, all $\mu_n$'s and $\nu_n$'s are
zero.

Finally, an application of \eqref{effe} to $x$ yields
$\tau_e\tau_f=1$ so that $\tau_e=\tau$ and $\tau_f=\tau^{-1}$ for
some $\tau\in\mathbb{C}\setminus\{0\}$.

We claim that all the actions corresponding to nonzero $\tau$ are
isomorphic to the specific action with $\tau=1$. The desired isomorphism is
given by the automorphism $\Phi_{\tau}:x\mapsto x$, $y\mapsto\tau y$. In
particular, $\left( \Phi_{\tau}\mathsf{e}_{\tau}\Phi_{\tau}^{-1}\right)
\left( y\right) =\tau^{-1}\Phi_{\tau}\left( \tau x\right)
=x=\mathsf{e}_{1}\left( y\right) $, where $\mathsf{e}_{\tau}\left( y\right)
$ denotes the action from (\ref{kkqq1}) with an arbitrary $\tau\neq0$.
\hfill\rule{0.5em}{0.5em}\smallskip

Now we consider the actions whose symbolic matrix $\left( \overset{\star
}{\mathsf{M}}_{\mathsf{ef}}\right) _{0}$ contains one~$\mathbf{\star}$.
Seemingly, \ the \ corresponding \ actions \ described \ below \ never \
appeared  in  the  literature before,  so we present more detailed
computations.\smallskip
\begin{theorem}
\label{Tb} The  $\left[ \left(
\begin{array}[c]{cc}
\mathbf{0} & \mathbf{\star}\\
\mathbf{0} & \mathbf{0}
\end{array}
\right)_0\!;\left(
\begin{array}[c]{cc}
\mathbf{0} & \mathbf{0}\\
\mathbf{0} & \mathbf{0}
\end{array}
\right)_1\right]$-series  consists  of  a  one-para\-meter
($b_0\in\mathbb{C}\setminus\{0\}$) family of $U_q(\mathfrak{sl}_2)$-module
algebra structures on the quantum plane
\begin{align}
\mathsf{k}(x) &=qx, & \mathsf{k}(y) &=q^{-2}y, &&&\label{kk2}\\
\mathsf{e}(x) &=0, & \mathsf{e}(y) &=b_0, &&&\label{kk3}\\
\mathsf{f}(x) &=b_0^{-1}xy, & \mathsf{f}(y) &=-qb_0^{-1}y^2. &&&\label{kk4}
\end{align}
All these structures are isomorphic, in particular to the action as above
with $b_0=1$.
\end{theorem}

The full action matrix of an action within this isomorphism class is of the
form
$$
\mathsf{M}=\left\Vert
\begin{array}[c]{cc}
qx & q^{-2}y\\
0 & 1\\
xy & -qy^2
\end{array}
\right\Vert .
$$\smallskip

P r o o f. First we demonstrate that an extension of
\eqref{kk2}--\eqref{kk4} to the entire action of $U_q(\mathfrak{sl}_2)$ on
$\mathbb{C}_q[x,y]$ passes through all the relations. It is clear that
\eqref{kk2} is compatible with the relation
$\mathsf{kk}^{-1}=\mathsf{k}^{-1}\mathsf{k}=\mathbf{1}$. Then we apply the
relations \eqref{ke}--\eqref{effe} to the quantum plane generators
\begin{align*}
(\mathsf{ke}-q^2\mathsf{ek})(x) &=\mathsf{k}(0)-q^3\mathsf{e}(x)=0,\\
(\mathsf{ke}-q^2\mathsf{ek})(y) &=\mathsf{k}(b_0)-\mathsf{e}(y)=b_0-b_{0}=0,
\\ (\mathsf{kf}-q^{-2}\mathsf{fk})(x)
&=\mathsf{k}\left(b_0^{-1}xy\right)-q^{-1}\mathsf{f}(x) \\
& =b_0^{-1}q^{-1}xy-q^{-1}b_0^{-1}xy=0, \\
(\mathsf{kf}-q^{-2}\mathsf{fk})(y)
&=\mathsf{k}\left(-qb_0^{-1}y^2\right)-q^{-4}\mathsf{f}(y) \\
&=-qb_0^{-1}q^{-4}y^2+q^{-4}\left(qb_0^{-1}y^2\right)=0, \\
\left(\mathsf{ef}-\mathsf{fe}-\dfrac{\mathsf{k}-
\mathsf{k}^{-1}}{q-q^{-1}}\right)(x)
&=\mathsf{e}\left(b_0^{-1}xy\right)-\mathsf{f}(0)-x
=b_0^{-1}\mathsf{e}(xy)-x \\
&=b_0^{-1}x\mathsf{e}(y)+b_0^{-1}\mathsf{e}(x)\mathsf{k}(y)-x=0, \\
\left(\mathsf{ef}-\mathsf{fe}-
\dfrac{\mathsf{k}-\mathsf{k}^{-1}}{q-q^{-1}}\right)(y)
&=-qb_0^{-1}\mathsf{e}\left(y^2\right)-\mathsf{f}(b_0)-
\dfrac{q^{-2}-q^2}{q-q^{-1}}y \\
&=-qb_0^{-1}\mathsf{e}\left(y^2\right)+\left(q+q^{-1}\right)y \\
&=-qb_0^{-1}y\mathsf{e}(y)-qb_0^{-1}\mathsf{e}(y)\mathsf{k}(y)+
\left(q+q^{-1}\right)y \\
&=-qy-q^{-1}y+\left(q+q^{-1}\right)y=0.
\end{align*}

Now apply the generators of $U_{2}\left( \mathfrak{sl}_{2}\right) $ to
(\ref{1}) and get%
\begin{align*}
\mathsf{k}\left(  yx-qxy\right)   &  =q^{-2}y\cdot qx-qqx\cdot q^{-2}y=0,\\
\mathsf{e}\left(  yx-qxy\right)   &  =y\mathsf{e}\left(  x\right)
+\mathsf{e}\left(  y\right)  \mathsf{k}\left(  x\right)  -qx\mathsf{e}\left(
y\right)  -q\mathsf{e}\left(  x\right)  \mathsf{k}\left(  y\right) \\
&  =0+b_{0}qx-qxb_{0}-0=0,\\
\mathsf{f}\left(  yx-qxy\right)   &  =\mathsf{f}\left(  y\right)
x+\mathsf{k}^{-1}\left(  y\right)  \mathsf{f}\left(  x\right)  -q\mathsf{f}%
\left(  x\right)  y-q\mathsf{k}^{-1}\left(  x\right)  \mathsf{f}\left(
y\right) \\
&  =-qb_{0}^{-1}y^{2}x+q^{2}yb_{0}^{-1}xy-qb_{0}^{-1}xy\cdot y+qq^{-1}x\cdot
qb_{0}^{-1}y^{2}\\
&  =-q^{3}b_{0}^{-1}xy^{2}+q^{3}b_{0}^{-1}xy^{2}-qb_{0}^{-1}xy^{2}+qb_{0}%
^{-1}xy^{2}=0.
\end{align*}

Next prove that $\left[\left(
\begin{array}[c]{cc}
\mathbf{0} & \mathbf{\star}\\
\mathbf{0} & \mathbf{0}
\end{array}
\right)_0;\left(
\begin{array}[c]{cc}
\mathbf{0} & \mathbf{0}\\
\mathbf{0} & \mathbf{0}
\end{array}
\right)_1\right]$-series contains no actions except
\eqref{kk2}--\eqref{kk4}. Show that the matrix elements of
$\mathsf{M}_{\mathsf{ef}}$ \eqref{mef} have no terms of degree higher than
two, viz. $(\mathsf{M}_{\mathsf{ef}})_n=0$ for $n\ge 3$. Now a general form
for $\mathsf{e}(x)$, $\mathsf{e}(y)$, $\mathsf{f}(x)$, $\mathsf{f}(y)$ is
\begin{align}
\mathsf{e}(x) &=\sum_{m+n\ge 0}\bar{\rho}_{mn}x^my^n, & \mathsf{e}(y)
&=\sum_{m+n\ge 0}\bar{\sigma}_{mn}x^my^n, & \label{er2}\\
\mathsf{f}(x) &=\sum_{m+n\ge 0}\bar{\rho}_{mn}'x^my^n, &
\mathsf{f}(y) &=\sum_{m+n\ge 0}\bar{\sigma}_{mn}'x^my^n & \label{fr2}
\end{align}
where $\bar{\rho}_{mn},\bar{\sigma}_{mn},\bar{\rho}_{mn}',
\bar{\sigma}_{mn}'\in\mathbb{C}$. Within this series one has the matrix of
weights
$$
\mathbf{wt}(\mathsf{M}_{\mathsf{ef}})=\left(
\begin{array}[c]{cc}
q^3 & 1\\
q^{-1} & q^{-4}
\end{array}
\right).
$$

In view of this, the general form \eqref{er2}--\eqref{fr2} should be a sum
of terms of the same weight
\begin{align}
\mathsf{e}(x) &=\sum_{m\ge 0}\rho_mx^{2m+3}y^m, \label{es2}\\
\mathsf{e}(y) &=b'+\sum_{m\ge 0}\sigma_mx^{2m+2}y^{m+1}, \label{es3}\\
\mathsf{f}(x) &=b''xy+\sum_{n\ge 0}\rho_n'x^{2n+3}y^{n+2}, \label{es4}\\
\mathsf{f}(y) &=b'''y^2+\sum_{n\ge 0}\sigma_n'x^{2n+2}y^{n+3}. \label{es5}
\end{align}

Now we combine \eqref{es2}--\eqref{es3}, \eqref{es4}--\eqref{es5}) with
\eqref{exy},  \eqref{fxy}, respectively, then project the resulting
relation to the one-dimensional subspace $\mathbb{C}x^{2m+3}y^{m+2}$ (resp.
$\mathbb{C}x^{2n+3}y^{n+3}$) (for every $m\ge 0$, resp. $n\ge 0$) to obtain
\begin{align*}
\dfrac{\rho_m}{\sigma_m} &=-q^2\dfrac{1-q^{m+1}}{1-q^{2m+4}} ,\\
\dfrac{\rho_n'}{\sigma_n'} &=-q^{-1}\dfrac{1-q^{n+3}}{1-q^{2n+4}} .
\end{align*}

Thus we get
\begin{align}
\mathsf{M}_{\mathsf{ef}} &=\left(
\begin{array}[c]{cc}
0 & b'\\
b^{''}xy & b^{'''}y^2
\end{array}
\right) \nonumber\\
& + \sum_{n\ge 0}\left(
\begin{array}[c]{cc}
\mu_{n}q^2(1-q^{n+1})x^{2n+3}y^n & -\mu_n(1-q^{2n+4})x^{2n+2}y^{n+1}\\
-\nu_n(1-q^{n+3})x^{2n+3}y^{n+2} & \nu_nq(1-q^{2n+4})x^{2n+2}y^{n+3}
\end{array}
\right) , \label{mp2}
\end{align}
where $\mu_n,\nu_n\in\mathbb{C}$. To prove that the second matrix vanishes,
assume the contrary. First consider the case when there exist both nonzero
$\mu_n$'s and $\nu_n$'s. As the sums here are finite, for the first row
choose the largest index $n_e$ with $\mu_{n_e}\ne 0$ and for the second
row, the largest index $n_f$ with $\nu_{n_f}\ne 0$. After applying
\eqref{definition}--\eqref{def1} one concludes that the highest degree of
monomials in $(ef-fe)(x)$ is $3n_e+3n_f+7$. This monomial is unique, and
its computation gives
\begin{equation}\label{mn}
\mu_{n_e}\nu_{n_f}q^{2n_en_f+2n_e}(1-q^{n_e+n_f+4})(1-q^{2n_e+2n_f+6})
x^{2n_e+2n_f+5}y^{n_e+n_f+2}.
\end{equation}

Under our assumptions on $q$, since $n_e\ge 0$, $n_f\ge 0$,
$\mu_{n_e}\nu_{n_f}\ne 0$, it becomes clear that \eqref{mn} is a nonzero
monomial of degree higher than $1$. This breaks \eqref{effe} whose r.h.s.
applied to $x$ has degree $1$. An application of \eqref{effe} to $x$ and
$y$ together with \eqref{mp2} leads to (up to terms of degree higher than
$1$)
\begin{align*}
\left(\mathsf{ef}-\mathsf{fe}-\dfrac{\mathsf{k}-
\mathsf{k}^{-1}}{q-q^{-1}}\right)(x) &=0=b'b''x-x,\\
\left(\mathsf{ef}-\mathsf{fe}-\dfrac{\mathsf{k}-
\mathsf{k}^{-1}}{q-q^{-1}}\right)(y)
&=0=b'b'''(1+q^{-2})y+\left(q+q^{-1}\right)y,
\end{align*}
which yields
$$b'=b_0,\qquad b''=b_0^{-1},\qquad b'''=-qb_0^{-1}$$
for some $b_0\ne 0$.

A similar, but simpler computation also shows that in the case
when all $\nu_n$'s are zero and some $\mu_n$'s are nonzero we have
the highest degree monomial of $(ef-fe)(x)$ of the form
$$
b_0^{-1}\mu_{n_e}\frac{(1-q^{n_e+3})(q^{2n_e+4}-1)}{1-q^2}
x^{2n_e+3}y^{n_e+1}.
$$
This monomial is nonzero due to our assumption on $q$, which gives the same
contradiction as above. The opposite case, when all $\mu_n$'s are zero and
some $\nu_n$'s are nonvanishing, can be treated similarly and also leads to
a contradiction. Therefore, all $\mu_n$'s and $\nu_n$'s are zero. This
gives the desired relations \eqref{kk2}--\eqref{kk4}.

Finally we show that the actions \eqref{kk2}--\eqref{kk4} with nonzero
$b_0$ are isomorphic to the specific action with $b_0=1$. The desired
isomorphism is as follows $\Phi_{b_0}:x\mapsto x$, $y\mapsto b_0y$. In
fact,
\begin{align*}
\left(\Phi_{b_0}\mathsf{e}_{b_0}\Phi_{b_0}^{-1}\right)(y)
&=\Phi_{b_0}\mathsf{e}_{b_0}\left(b_0^{-1}y\right)=b_0^{-1}\Phi_{b_0}(b_0)
=\Phi_{b_0}(1)=1=\mathsf{e}_1(y), & \\
\left(\Phi_{b_0}\mathsf{f}_{b_0}\Phi_{b_0}^{-1}\right)(x)
&=\Phi_{b_0}\mathsf{f}_{b_0}(x)=b_0^{-1}\Phi_{b_{0}}(xy)=b_0^{-1}b_0xy=xy
=\mathsf{f}_1(x) , & \\
\left(\Phi_{b_0}\mathsf{f}_{b_0}\Phi_{b_0}^{-1}\right)(y)
&=\Phi_{b_0}\mathsf{f}_{b_0}\left(b_0^{-1}y\right)
=b_0^{-1}\Phi_{b_0}\left(-qb_0^{-1}y^2\right)=-qb_0^{-2}b_0^2y^2=
& \\
&=-qy^2=\mathsf{f}_1(y) . &
\end{align*}
The theorem is proved. \hfill\rule{0.5em}{0.5em}

\begin{theorem}
The $\left[\left(
\begin{array}
[c]{cc}
\mathbf{0} & \mathbf{0}\\
\mathbf{\star} & \mathbf{0}
\end{array}
\right)_0;\left(
\begin{array}
[c]{cc}
\mathbf{0} & \mathbf{0}\\
\mathbf{0} & \mathbf{0}
\end{array}
\right)_1\right]$-series consists of a one-para\-meter
($c_0\in\mathbb{C}\setminus\{0\}$) family of $U_q(\mathfrak{sl}_2)$-module
algebra structures on the quantum plane
\begin{align}
\mathsf{k}(x) &=q^2x, & \mathsf{k}(y) &=q^{-1}y, &&&\label{kc}\\
\mathsf{e}(x) &=-qc_0^{-1}x^2, & \mathsf{e}(y) &=c_0^{-1}xy, &&&\label{ec}\\
\mathsf{f}(x) &=c_0, & \mathsf{f}(y) &=0. &&&\label{fc}
\end{align}

All these structures are isomorphic, in particular to the action as above
with $c_0=1$.
\end{theorem}

The full action matrix for this isomorphism class (with $c_0=1$) is
$$
\mathsf{M}=\left\Vert
\begin{array}[c]{cc}
q^2x & q^{-1}y\\
-qx^2 & xy\\
1 & 0
\end{array}
\right\Vert .
$$

P r o o f. \ Quite literally repeats that of the previous theorem.
\hfill\rule{0.5em}{0.5em}
\begin{theorem}
\label{T3h}The $\left[\left(
\begin{array}
[c]{cc}
\mathbf{\star} & \mathbf{0}\\
\mathbf{0} & \mathbf{0}
\end{array}
\right)_0;\left(
\begin{array}
[c]{cc}
\mathbf{0} & \mathbf{0}\\
\mathbf{0} & \mathbf{0}
\end{array}
\right)_1\right]$-series consists of a three-para\-meter
($a_0\in\mathbb{C}\setminus\{0\}$, $s,t\in\mathbb{C}$) family of
$U_q(\mathfrak{sl}_2)$-actions on the quantum plane
\begin{align}
\mathsf{k}(x) &=q^{-2}x, & \mathsf{k}(y) &=q^{-1}y, &&&\label{k11}\\
\mathsf{e}(x) &=a_0, & \mathsf{e}(y) &=0, &&&\label{k12}\\
\mathsf{f}(x) &=-qa_0^{-1}x^2+ty^4, & \mathsf{f}(y) &=-qa_0^{-1}xy+sy^3.
&&&\label{k13}
\end{align}

The generic domain $\{(a_0,s,t)|\:s\ne 0,\;t\ne 0\}$ with respect
to the parameters splits into uncountably many disjoint subsets
$\{(a_0,s,t)|\:s\ne 0,\;t\ne 0,\;\varphi=\mathrm{const}\}$, where
$\varphi=\dfrac{t}{a_0s^2}$. Each of those subsets corresponds to
an isomorphism class of $U_q(\mathfrak{sl}_2)$-module algebra
structures. Additionally, there exist three more isomorphism
classes corresponding to the subsets
$$
\{(a_0,s,t)|\:s\ne 0,t=0\},\quad\{(a_0,s,t)|s=0,t\ne 0\},\quad
\{(a_0,s,t)|\:s=0,t=0\}.
$$
\end{theorem}

P r o o f. \ A routine verification demonstrates that
\eqref{k11}--\eqref{k13} pass through all the relations as before, hence
admit an extension to a well-defined series of
$U_q(\mathfrak{sl}_2)$-actions on the quantum plane.\smallskip

Now check that $\left[\left(
\begin{array}[c]{cc}
\mathbf{\star} & \mathbf{0}\\
\mathbf{0} & \mathbf{0}
\end{array}
\right)_0;\left(
\begin{array}[c]{cc}
\mathbf{0} & \mathbf{0}\\
\mathbf{0} & \mathbf{0}
\end{array}
\right)_1\right]$-series contains no other actions except
\eqref{k11}--\eqref{k13}. First consider the polynomial $\mathsf{e}(x)$.
Since its weight is $q^2\mathbf{wt}(x)=1$, and the weight of any monomial
other than constant is a negative degree of $q$ (within the series under
consideration), hence not $1$, one gets $\mathsf{e}(x)=a_0$. In a similar
way, the only possibility for $\mathsf{e}(y)$ is zero, because if not,
$\mathbf{wt}(\mathsf{e}(y))=q^2\mathbf{wt}(y)=q$, which is impossible in
view of the above observations.

Turn to $\mathsf{f}\left( x\right) $ and observe that $\mathbf{wt}\left(
\mathsf{f}\left( x\right) \right) =q^{-4}$. It is easy to see that all the
monomials with this weight are $x^{2}$, $xy^{2}$, $y^{4}$, that is
$\mathsf{f}\left( x\right) =ux^{2}+vxy^{2}+wy^{4}$. In a similar way
$\mathbf{wt}\left( \mathsf{f}\left( y\right) \right) =q^{-3}$ and so
$\mathsf{f}\left( y\right) =zxy+sy^{3}$. A substitution to (\ref{effe})
yields $\left( 1+q^{-2}\right) ua_{0}=-\left( q+q^{-1}\right) $, $v=0$,
$za_{0}q^{-1}=-1$. Note that (\ref{fxy}) gives no new relations for $u$,
$v$, $z$ and provides no restriction on $w$ and $s$ at all. This leads to
(\ref{k13}).

To distinguish the isomorphism classes of the structures within
this series, we use Theorem \ref{T1} in writing down the general
form of an automorphism of $\mathbb{C}_q[x,y]$ as
$\Phi_{\theta,\omega}:x\mapsto\theta x$, $y\mapsto\omega y$.
Certainly, this commutes with the action of $\mathsf{k}$. For
other generators we get
\begin{align*}
\left(\Phi_{\theta,\omega}\mathsf{e}_{a_0,s,t}\Phi_{\theta,\omega}^{-1}
\right)(x) &=\Phi_{\theta,\omega}\mathsf{e}_{a_0,s,t}\left(\theta^{-1}x
\right)=\theta^{-1}a_0, \\
\left(\Phi_{\theta,\omega}\mathsf{e}_{a_0,s,t}\Phi_{\theta,\omega}^{-1}
\right)(y) &=\Phi_{\theta,\omega}\mathsf{e}_{a_0,s,t}\left(\omega^{-1}y
\right)=\omega^{-1}\Phi_{\theta,\omega}\mathsf{e}_{a_0,s,t}(y)=0, \\
\left(\Phi_{\theta,\omega}\mathsf{f}_{a_0,s,t}\Phi_{\theta,\omega}^{-1}
\right)(x) &=\Phi_{\theta,\omega}\mathsf{f}_{a_0,s,t}
\left(\theta^{-1}x\right)=\theta^{-1}\Phi_{\theta,\omega}
\left(-qa_0^{-1}x^2+ty^4\right) \\
&=-qa_0^{-1}\theta x^2+\theta^{-1}t\omega^4y^4, \\
\left(\Phi_{\theta,\omega}\mathsf{f}_{a_0,s,t}\Phi_{\theta,\omega}^{-1}
\right)(y) &=\Phi_{\theta,\omega}\mathsf{f}_{a_0,s,t}\left(\omega^{-1}y
\right)=\omega^{-1}\Phi_{\theta,\omega}\left(-qa_0^{-1}xy+sy^3\right) \\
&=-q\theta a_0^{-1}xy+s\omega^2y^3.
\end{align*}
That is, the automorphism $\Phi_{\theta,\omega}$ transforms the parameters
of actions \eqref{k12}--\eqref{k13} as follows:
$$
a_0\mapsto\theta^{-1}a_0,\qquad s\mapsto\omega^2s,\qquad
t\mapsto\theta^{-1}\omega^4t.
$$
In particular, this means that within the domain $\{s\ne 0,\;t\ne 0\}$ one
obtains an invariant $\varphi=\dfrac{t}{a_0s^2}$ of the isomorphism class.
Obviously, the complement to this domain further splits into three distinct
subsets $\{s\ne 0,\;t=0\}$, $\{s=0,\;t\ne 0\}$, $\{s=0,\;t=0\}$
corresponding to the isomorphism classes listed in the formulation, and our
result follows. \hfill\rule{0.5em}{0.5em}
\medskip

Note that up to isomorphism of $U_q(\mathfrak{sl}_2)$-module algebra
structure, the full action matrix corresponding to \eqref{k11}--\eqref{k13}
is of the form
$$
\mathsf{M}=\left\Vert
\begin{array}[c]{cc}
q^{-2}x & q^{-1}y\\
1 & 0\\
-qx^2+ty^4 & -qxy+sy^3
\end{array}
\right\Vert .
$$\smallskip
\begin{theorem}
\label{T3l}The $\left[\left(
\begin{array}
[c]{cc}
\mathbf{0} & \mathbf{0}\\
\mathbf{0} & \mathbf{\star}
\end{array}
\right)_0;\left(
\begin{array}
[c]{cc}
\mathbf{0} & \mathbf{0}\\
\mathbf{0} & \mathbf{0}
\end{array}
\right)_1\right]$-series consists of three-parameter
($d_0\in\mathbb{C}\setminus\{0\}$, $s,t\in\mathbb{C}$) family of
$U_q(\mathfrak{sl}_2)$-actions on the quantum plane
\begin{align}
\mathsf{k}(x) &=qx, & \mathsf{k}(y) &=q^2y, &&&\label{k31}\\
\mathsf{e}(x) &=-qd_0^{-1}xy+sx^3, & \mathsf{e}(y) &=-qd_0^{-1}y^2+tx^4,
&&&\label{k32}\\
\mathsf{f}(x) &=0, & \mathsf{f}(y) &=d_0. &&&\label{k33}
\end{align}

Here we have the domain $\{(d_0,s,t)|\:s\ne 0,\;t\ne 0\}$ which
splits into the disjoint subsets $\{(d_0,s,t)|\:s\ne 0,\;t\ne
0,\;\varphi=\mathrm{const}\}$ with $\varphi=\dfrac{t}{d_0s^2}$.
This uncountable family of subsets is in one-to-one correspondence
to the isomorphism classes of $U_q(\mathfrak{sl}_2)$-module
algebra structures. Aside of those, one also has three more
isomorphism classes labelled by the subsets $\{(d_0,s,t)|\:s\ne
0,\;t=0\}$, $\{(d_0,s,t)|\:s=0,\;t\ne 0\}$,
$\{(d_0,s,t)|\:s=0,\;t=0\}$.
\end{theorem}

P r o o f. \ Is the same as that of the previous theorem.
\hfill\rule{0.5em}{0.5em}\medskip 

Here, also up to isomorphism of $U_q(\mathfrak{sl}_2)$-module
algebra structures, the full action matrix is
$$
\mathsf{M}=\left\Vert
\begin{array}[c]{cc}
qx & q^2y\\
-qxy+sx^3 & -qy^2+tx^4\\
0 & 1
\end{array}
\right\Vert .
$$

\begin{remark}
{\rm There could be no isomorphisms between the $U_{q}\left(  \mathfrak{sl}%
_{2}\right) $-module algebra structures on $\mathbb{C}_{q}[x,y]$ picked
from different series. This is because every automorphism of the quantum
plane commutes with the action of $\mathsf{k}$, hence, the restrictions of
isomorphic actions to $\mathsf{k}$ are always the same. On the other hand,
the actions of $\mathsf{k}$ in different series are different.}
\end{remark}
\begin{remark}
{\rm The list of $U_{q}\left( \mathfrak{sl}_{2}\right) $-module algebra
structures on $\mathbb{C}_{q}[x,y]$ presented in the theorems of this
section is complete. This is because the assumptions of those theorems
exhaust all admissible forms for the components $\left(
\mathsf{M}_{\mathsf{ef}}\right) _{0}$, $\left(
\mathsf{M}_{\mathsf{ef}}\right) _{1}$ of the action
$\mathsf{ef}$-matrix.\textit{\ }}
\end{remark}
\begin{remark}
{\rm In all series of $U_q(\mathfrak{sl}_2)$-module algebra structures
listed in Theorems \ref{T01}--\ref{T3l}, except the series $\left[\left(
\begin{array}[c]{cc}
\mathbf{0} & \mathbf{0}\\
\mathbf{0} & \mathbf{0}
\end{array}
\right)_0;\left(
\begin{array}[c]{cc}
\mathbf{0} & \mathbf{0}\\
\mathbf{0} & \mathbf{0}
\end{array}
\right)_1\right]$, the weight constants $\alpha$ and $\beta$ satisfy the
assumptions of Proposition \ref{Pmono}. So the claim of this proposition is
well visible in a rather simple structure of nonzero homogeneous components
of $\mathsf{e}(x)$, $\mathsf{e}(y)$, $\mathsf{f}(x)$, $\mathsf{f}(y)$,
which everywhere reduce to monomials.}
\end{remark}

{\center\section{\label{Comp}Composition Series}}

Let us view the $U_{q}\left( \mathfrak{sl}_{2}\right) $-module
algebra structures on $\mathbb{C}_{q}[x,y]$ listed in the theorems
of the previous section merely as representations of $U_{q}\left(
\mathfrak{sl}_{2}\right) $ in the \textit{vector space}
$\mathbb{C}_{q}[x,y]$. Our immediate intention is to describe the
composition series for these representations.
\begin{proposition}
The representations corresponding to $\left[\!\left(\!\!
\begin{array}[c]{cc}
\mathbf{0} & \mathbf{0}\\
\mathbf{0} & \mathbf{0}
\end{array}
\!\!\right)_{\!\!0}\!\!;\left(\!\!
\begin{array}[c]{cc}
\mathbf{0} & \mathbf{0}\\
\mathbf{0} & \mathbf{0}
\end{array}
\!\!\right)_{\!\!1}\right]$-series described in \eqref{t1}--\eqref{t2}
split into the direct sum
$\mathbb{C}_q[x,y]=\oplus_{m=0}^\infty\oplus_{n=0}^\infty\mathbb{C}x^my^n$
of (irreducible) one-dimensional subrepresentations. These
subrepresentations may belong to two isomorphism classes, depending on the
weight of a specific monomial $x^my^n$ which can be $\pm 1$ (see
Th.~\ref{T01}).
\end{proposition}
P r o o f. \ Since $\mathsf{e}$ and $\mathsf{f}$ are represented by zero
operators and the monomials $x^my^n$ are eigenvectors for $\mathsf{k}$,
then every direct summand is $U_q(\mathfrak{sl}_2)$-invariant.
\hfill\rule{0.5em}{0.5em}\smallskip

Now turn to nontrivial $U_{q}\left( \mathfrak{sl}_{2}\right) $-module
algebra structures and start with the well-known case \cite{kassel,
mon/smi}.
\begin{proposition}
The representations corresponding to $\left[\!\left(\!\!
\begin{array}[c]{cc}
\mathbf{0} & \mathbf{0}\\
\mathbf{0} & \mathbf{0}
\end{array}
\!\!\right)_{\!\!0}\!\!;\left(\!\!
\begin{array}[c]{cc}
\mathbf{0} & \mathbf{\star}\\
\mathbf{\star} & \mathbf{0}
\end{array}
\!\!\right)_{\!\!1}\right]$-series described in \eqref{kkqq}--\eqref{ffxy}
split into the direct sum
$\mathbb{C}_q[x,y]=\oplus_{n=0}^\infty\mathbb{C}_q[x,y]_n$ of irreducible
finite-dimensional subrepresentations, where $\mathbb{C}_q[x,y]_n$ is the
$n$-th homogeneous component (introduced in Sect.~3) 
with $\dim\mathbb{C}_q[x,y]_n=n+1$ and the isomorphism class of this
subrepresentation is $\mathcal{V}_{1,n}$ \cite[Ch.~VI]{kassel}.
\end{proposition}

P r o o f. \ Is that of Theorem VII.3.3 (b) from \cite{kassel}.
\hfill\rule{0.5em}{0.5em}\smallskip

In the subsequent observations we encounter a split picture which does not
reduce to a collection of purely finite-dimensional sub- or quotient
modules. We recall the definition of the Verma modules in our specific case
of $U_{q}\left( \mathfrak{sl}_{2}\right) $.
\begin{definition}
A Verma module $\mathcal{V}(\lambda)$
($\lambda\in\mathbb{C}\setminus\{0\}$) is a vector space with a~basis
$\{v_i,\:i\ge 0\}$, where the $U_q(\mathfrak{sl}_2)$ action is given by
\begin{align*}
\mathsf{k}v_i &=\lambda q^{-2i}v_{i},\qquad
\mathsf{k}^{-1}v_i=\lambda^{-1}q^{2i}v_{i}, & \\
\mathsf{e}v_0 &=0,\qquad
\mathsf{e}v_{i+1}=\dfrac{\lambda q^{-i}-\lambda^{-1}q^i}{q-q^{-1}}v_i,
\qquad\mathsf{f}v_i=\dfrac{q^{i+1}-q^{-i-1}}{q-q^{-1}}v_{i+1}. &
\end{align*}
\end{definition}

Note that the Verma module $\mathcal{V}\left( \lambda\right) $ is generated
by the highest weight vector $v_{0}$ whose weight is $\lambda$ (for details
see, e.g., \cite{kassel}).
\begin{proposition}\label{Ph}
The representations corresponding to $\left[\!\left(\!\!
\begin{array}[c]{cc}
\mathbf{0} & \mathbf{\star}\\
\mathbf{0} & \mathbf{0}
\end{array}
\!\!\right)_{\!\!0}\!\!;\left(\!\!
\begin{array}[c]{cc}
\mathbf{0} & \mathbf{0}\\
\mathbf{0} & \mathbf{0}
\end{array}
\!\!\right)_{\!\!1}\right]$-series described in \eqref{kk1}--\eqref{kk4}
split into the direct sum of subrepresentations
$\mathbb{C}_q[x,y]=\oplus_{n=0}^\infty\mathcal{V}_n$, where
$\mathcal{V}_n=x^n\mathbb{C}[y]$. Each $\mathcal{V}_n$ admits a composition
series of the form $0\subset\mathcal{J}_n\subset\mathcal{V}_n$. The simple
submodule $\mathcal{J}_n$ of dimension $n+1$ is the linear span of
$x^n,x^ny,\ldots,x^ny^{n-1},x^ny^n$, whose isomorphism class is
$\mathcal{V}_{1,n}$ and $\mathcal{J}_n$ is not a direct summand in the
category of $U_q(\mathfrak{sl}_2)$-modules (there exist no submodule
$\mathcal{W}$ such that $\mathcal{V}_n=\mathcal{J}_n\oplus\mathcal{W}$).
The quotient module $\mathcal{V}_n\diagup \mathcal{J}_n=\mathcal{Z}_n$ is
isomorphic to the (simple) Verma module $\mathcal{V}\left(q^{-n-2}\right)$.
\end{proposition}

P r o o f. \ Due to the isomorphism statement of Theorem \ref{Tb}, it
suffices to set the parameter of the series $b_0=1$ in
\eqref{kk1}--\eqref{kk4}. An application of $\mathsf{e}$ and $\mathsf{f}$
to the basis elements of $\mathbb{C}_q[x,y]$ gives
\begin{align}
\mathsf{e}(x^ny^p) &=q^{1-p}\dfrac{q^{p}-q^{-p}}{q-q^{-1}}x^ny^{p-1}\ne 0,
\qquad\forall p>0, \label{exyp}\\
\mathsf{e}(x^n) &=0, \label{exn}\\
\mathsf{f}(x^ny^p) &=q^{-n}\dfrac{q^{2n}-q^{2p}}{q-q^{-1}}x^ny^{p+1},\qquad
\forall p\ge 0, \label{fxyp}
\end{align}
which already implies that each $\mathcal{V}_n$ is
$U_q(\mathfrak{sl}_2)$-invariant. Also $\mathcal{J}_n$ is a submodule of
$\mathcal{V}_n$ generated by the highest weight vector $x^n$, as the
sequence of weight vectors $\mathsf{f}(x^ny^p)$ terminates because
$\mathsf{f}(x^ny^n)=0$. The highest weight of $\mathcal{J}_n$ is $q^n$,
hence by Theorem VI.3.5 of \cite{kassel}, the submodule $\mathcal{J}_n$ is
simple and its isomorphism class is $\mathcal{V}_{1,n}$.

Now assume the contrary to our claim, that is
$\mathcal{V}_n=\mathcal{J}_n\oplus\mathcal{W}$ for some submodule
$\mathcal{W}$ of $\mathcal{V}_n$, and $\mathcal{V}_n\ni x^ny^{n+1}=u+w$,
$u\in\mathcal{J}_n$, $w\in\mathcal{W}$ is the associated decomposition. In
view of \eqref{exyp}--\eqref{exn}, an application of $\mathsf{e}^{n+1}$
gives $A(q)x^n=\mathsf{e}^{n+1}(w)$ for some nonzero constant $A(q)$,
because $\mathsf{e}^{n+1}|_{\mathcal{J}_n}=0$. This is a contradiction,
because $\mathcal{J}_n\cap\mathcal{W}=\{0\}$, thus there exist no submodule
$\mathcal{W}$ as above.

The quotient module $\mathcal{Z}_n$ is spanned by its basis vectors
$z_{n+1,}z_{n+2},\ldots$ which are the projections of
$x^ny^{n+1},x^ny^{n+2},\ldots$ respectively, to
$\mathcal{V}_n\diagup\mathcal{J}_n$. It follows from \eqref{exyp}, that
$z_{n+1}$ is the highest weight vector whose weight is $q^{-n-2}$, and it
generates $\mathcal{Z}_n$ by \eqref{fxyp}. Now the universality property of
the Verma modules (see, e.g., \cite[Prop.~VI.3.7]{kassel}) implies that
there exists a surjective morphism of modules $\Pi:\mathcal{V}\left(
q^{-n-2}\right)\to\mathcal{Z}_n$. It follows from Proposition 2.5 of
\cite{jantzen} that $\ker\Pi=0$, hence $\Pi$ is an isomorphism.
\hfill\rule{0.5em}{0.5em}\smallskip

The next series, unlike the previous one, involves the lowest
weight Verma modules. In all other respects the proof of the
following proposition is the same (we also set here $d_{0}=1$).

\begin{proposition}
The representations corresponding to $\left[\!\left(\!\!
\begin{array}[c]{cc}
\mathbf{0} & \mathbf{0}\\
\mathbf{\star} & \mathbf{0}
\end{array}
\!\!\right)_{\!\!0}\!\!;\left(\!\!
\begin{array}[c]{cc}
\mathbf{0} & \mathbf{0}\\
\mathbf{0} & \mathbf{0}
\end{array}
\!\!\right)_{\!\!1}\right]$-series described in \eqref{kc}--\eqref{fc}
split into the direct sum of subrepresentations
$\mathbb{C}_q[x,y]=\oplus_{n=0}^\infty\mathcal{V}_n$, where
$\mathcal{V}_n=\mathbb{C}[x]y^n$. Each $\mathcal{V}_n$ admits a composition
series of the form $0\subset\mathcal{J}_n\subset\mathcal{V}_n$. The simple
submodule $\mathcal{J}_n$ of dimension $n+1$ is the linear span of
$y^n,xy^n,\ldots,x^{n-1}y^n,x^ny^n$. This is a finite-dimensional
$U_q(\mathfrak{sl}_2)$-module whose lowest weight vector is $y^n$ with
weight $q^{-n}$, and its isomorphism class is $\mathcal{V}_{1,n} $. Now the
submodule $\mathcal{J}_n$ is not a direct summand in the category of
$U_q(\mathfrak{sl}_2)$-modules (there exists no submodule $\mathcal{W}$
such that $\mathcal{V}_n=\mathcal{J}_n\oplus\mathcal{W}$). The quotient
module $\mathcal{V}_n\diagup\mathcal{J}_n=\mathcal{Z}_n$ is isomorphic to
the (simple) Verma module with lowest weight $q^{n+2}$.
\end{proposition}

Now turn to considering the three parameter series as in Theorems
\ref{T3h}, \ref{T3l}. Despite we have now three parameters, the
entire series has the same split picture.

\begin{proposition}
The representations corresponding to $\left[\!\left(\!\!
\begin{array}[c]{cc}
\mathbf{\star} & \mathbf{0}\\
\mathbf{0} & \mathbf{0}
\end{array}
\!\!\right)_{\!\!0}\!\!;\left(\!\!
\begin{array}[c]{cc}
\mathbf{0} & \mathbf{0}\\
\mathbf{0} & \mathbf{0}
\end{array}
\!\!\right)_{\!\!1}\right]$-series described in \eqref{k11}--\eqref{k13}
split into the direct sum of subrepresentations
$\mathbb{C}_q[x,y]=\oplus_{n=0}^\infty\mathcal{V}_n$, where $\mathcal{V}_n$
is a submodule generated by its highest weight vector $y^n$. Each
$\mathcal{V}_n$ with $n\ge 1$ is isomorphic to a simple highest weight
Verma module $\mathcal{V}\left(q^{-n}\right)$. The submodule
$\mathcal{V}_0$ admits a composition series of the form
$0\subset\mathcal{J}_0\subset\mathcal{V}_0$, where
$\mathcal{J}_0=\mathbb{C}\mathbf{1}$. The submodule $\mathcal{J}_0$ is not
a direct summand in the category of $U_q(\mathfrak{sl}_2)$-modules (there
exists no submodule $\mathcal{W}$ such that
$\mathcal{V}_0=\mathcal{J}_0\oplus\mathcal{W}$). The quotient module
$\mathcal{V}_0\diagup\mathcal{J}_0$ is isomorphic to the (simple) Verma
module $\mathcal{V}\left(q^{-2}\right)$.
\end{proposition}

P r o o f. \ First, let us consider the special case of \eqref{k12},
\eqref{k13} in which $s=t=0$ and $a_{0}=1$. Then
$\mathcal{V}_n=\mathbb{C}[x]y^n$ are $U_q(\mathfrak{sl}_2)$-invariant, and
we calculate
\begin{align}
\mathsf{e}(x^py^n) &=q^{-n-p+1}\dfrac{q^p-q^{-p}}{q-q^{-1}}x^{p-1}y^n\ne 0,
\qquad\forall p>0, & \nonumber \\
\mathsf{e}(y^n) &=0, & \nonumber \\
\mathsf{f}(x^py^n) &=q^{n+p}\dfrac{q^{p+n}-q^{-p-n}}{q-q^{-1}}x^{p+1}y^n,
\qquad\forall p\ge 0. & \label{fxy1}
\end{align}

Note that $\mathsf{f}(x^py^n)=0$ only when $p=n=0$. Therefore
$\mathcal{V}_n$ admits a generating highest weight vector $y^n$
whose weight is $q^{-n}$. As in the proof of Proposition \ref{Ph}
we deduce that each $\mathcal{V}_n$ with $n\ge 1$ is isomorphic to
the (highest weight simple) Verma module
$\mathcal{V}\left(q^{-n}\right)$. In the case $n=0$, it is clear
that $\mathcal{V}_0$ contains an obvious submodule
$\mathbb{C}\mathbf{1}$ which is not a direct summand by an
argument in the proof of Proposition \ref{Ph}.

Turn to the general case when the three parameters are unrestricted. The
formulas \eqref{k11}--\eqref{k13} imply the existence of a descending
sequence of submodules
$$
\ldots\subset\mathcal{F}_{n+1}\subset\mathcal{F}_{n}\subset
\mathcal{F}_{n-1}\subset\ldots\subset\mathcal{F}_2\subset\mathcal{F}_1
\subset\mathcal{F}_0=\mathbb{C}_q[x,y],
$$
where $\mathcal{F}_n=\cup_{k=n}^\infty\mathbb{C}[x]y^k$, because operators
of the action, being applied to a monomial, can only increase its degree in
$y$. Note that the quotient module $\mathcal{F}_n\diagup\mathcal{F}_{n+1}$
with unrestricted parameters is isomorphic to the module
$\mathbb{C}[x]y^n\cong\mathcal{V}\left(q^{-n}\right)$, just as in the
case $s=t=0$.

Now we claim that $\mathcal{F}_{n+1}$ is a direct summand in
$\mathcal{F}_{n} $, namely
$\mathcal{F}_{n}=\mathcal{V}_{n}\oplus\mathcal{F}_{n+1}$, $n\geq0 $, with
$\mathcal{V}_{n}=U_{q}\left( \mathfrak{sl}_{2}\right) y^{n}$ for $n\geq1$
and $\mathcal{V}_{0}=U_{q}\left( \mathfrak{sl}_{2}\right) x$.

First consider the case $n\ge 1$. By virtue of \eqref{k11}--\eqref{k13},
$y^n$ is a generating highest weight vector of the submodule
$\mathcal{V}_n=U_q(\mathfrak{sl}_2)y^n$, whose weight is $q^{-n}$. Another
application of the argument in the proof of Proposition \ref{Ph}
establishes an isomorphism
$\mathcal{V}_n\cong\mathcal{V}\left(q^{-n}\right)$; in particular,
$\mathcal{V}_n$ is a simple module by Proposition 2.5 of \cite{jantzen}.
Hence $\mathcal{V}_n\cap\mathcal{F}_{n+1}$ can not be a proper submodule of
$\mathcal{V}_n$. Since $\mathcal{V}_n$ is not contained in
$\mathcal{F}_{n+1}$ (as $y^n\notin\mathcal{F}_{n+1}$), the latter
intersection is zero, and the sum $\mathcal{V}_n+\mathcal{F}_{n+1}$ is
direct. On the other hand, a comparison of \eqref{k13} and \eqref{fxy1}
allows one to deduce that $\mathcal{V}_n+\mathcal{F}_{n+1}$ contains all
the monomials $x^py^m$, $m\ge n$, $p\ge 0$. This already proves
$\mathcal{F}_n=\mathcal{V}_n\oplus\mathcal{F}_{n+1}$.

Turn to the case $n=0$. The composition series
$0\subset\mathbb{C}\mathbf{1}\subset\mathcal{V}_0=U_q(\mathfrak{sl}_2)x$
is treated in the same way as that for $\mathcal{V}_0$ in
Proposition \ref{Ph}; in particular, the quotient module
$\mathcal{V}_0/\mathbb{C}\mathbf{1}$ is isomorphic to the simple
Verma module $\mathcal{V}\left(q^{-2}\right)$. Let
$\pi:\mathcal{V}_0\to\mathcal{V}_0/\mathbb{C}\mathbf{1}$ be the
natural projection map. Obviously, $\mathcal{F}_1$ does not
contain $\mathbb{C}\mathbf{1}$, hence the restriction of $\pi$ to
$\mathcal{V}_0\cap\mathcal{F}_1$ is one-to-one. Thus, to prove
that the latter intersection is zero, it suffices to verify that
$\pi(\mathcal{V}_0\cap\mathcal{F}_1)$ is zero. As the module
$\mathcal{V}_0/\mathbb{C}\mathbf{1}$ is simple, the only
alternative to $\pi(\mathcal{V}_0\cap\mathcal{F}_1)=\{0\}$ could
be
$\pi(\mathcal{V}_0\cap\mathcal{F}_1)=\mathcal{V}_0/\mathbb{C}\mathbf{1}$.
Under the latter assumption, there should exist some element of
$\mathcal{V}_0\cap\mathcal{F}_1$, which is certainly of the form
$Py$ for some $P\in\mathbb{C}_q[x,y]$, and such that
$\pi(x)=\pi(Py)$. This relation is equivalent to $x-Py=\gamma$ for
some constant $\gamma$, which is impossible, because the monomials
that form $Py$, together with $x$ and $\mathbf{1}$, are linearly
independent. The contradiction we get this way proves that
$\mathcal{V}_0\cap\mathcal{F}_1=\{0\}$, hence the sum
$\mathcal{V}_0+\mathcal{F}_1$ is direct. On the other hand, a
comparison of (\ref{k13}) and (\ref{fxy1}) allows one to deduce
that $\mathcal{V}_0+\mathcal{F}_1$ contains all the monomials
$x^py^m$, with $m,p\ge 0$. Thus the relation
$\mathcal{F}_n=\mathcal{V}_n\oplus\mathcal{F}_{n+1}$ is now proved
for all $n\ge 0$. This, together with
$\cap_{i=0}^\infty\mathcal{F}_i=\{0\}$, implies that
$$
\mathbb{C}_q[x,y]=\left(\oplus_{n=1}^\infty U_q(\mathfrak{sl}_2)y^n\right)
\oplus U_q(\mathfrak{sl}_2)x,
$$
which was to be proved. \hfill\rule{0.5em}{0.5em}\smallskip

In a similar way we obtain the following

\begin{proposition}
The representations corresponding to $\left[\!\left(\!\!
\begin{array}[c]{cc}
\mathbf{0} & \mathbf{0}\\
\mathbf{0} & \mathbf{\star}
\end{array}
\!\!\right)_{\!\!0};\left(\!\!
\begin{array}[c]{cc}
\mathbf{0} & \mathbf{0}\\
\mathbf{0} & \mathbf{0}
\end{array}
\!\!\right)_{\!\!1}\right]$-series described in \eqref{k31}--\eqref{k33}
split into the direct sum of subrepresentations
$\mathbb{C}_q[x,y]=\oplus_{n=0}^\infty\mathcal{V}_n$, where $\mathcal{V}_n$
is a submodule generated by its lowest weight vector $x^n$. Each
$\mathcal{V}_n$ with $n\ge 1$ is isomorphic to a simple lowest weight Verma
module whose lowest weight is $q^{n}$. The submodule $\mathcal{V}_0$ admits
a composition series of the form
$0\subset\mathcal{J}_0\subset\mathcal{V}_0$, where
$\mathcal{J}_0=\mathbb{C}\mathbf{1}$. The submodule $\mathcal{J}_0$ is not
a direct summand in the category of $U_q(\mathfrak{sl}_2)$-modules (there
exists no submodule $\mathcal{W}$ such that
$\mathcal{V}_0=\mathcal{J}_0\oplus\mathcal{W}$). The quotient module
$\mathcal{V}_0\diagup\mathcal{J}_0$ is isomorphic to the (simple) lowest
weight Verma module whose lowest weight is $q^2$.
\end{proposition}
\indent The associated classical limit actions of the Lie algebra
$\mathfrak{sl}_2$ (here it is the Lie algebra generated by $e$, $f$, $h$
subject to the relations $[h,e]=2e$, $[h,f]=-2f$, $[e,f]=h$) on
$\mathbb{C}[x,y]$ by differentiations is derived from the quantum action
via substituting $k=q^h$ with subsequent formal passage to the limit as
$q\to 1$.

In this way we present all quantum and classical actions in Table 1. It
should be noted  that  there  exist  more $\mathfrak{sl}_2$-actions on
$\mathbb{C}[x,y]$ by differentiations (see, e.g., \cite{GLKO}) than one can
see in Table 1. It follows from our results that the rest of the classical
actions admit no quantum counterparts. On the other hand, among the quantum
actions listed in the first row of Table 1, the only one to which the above
classical limit procedure is applicable, is the action with $k(x)=x$,
$k(y)=y$. \newline The rest three actions of this series admit no classical
limit in
the above sense.\medskip

\emph{Acknowledgements}. One of the authors (S.D.) is thankful to Yu.
Bespalov, J. Cuntz, B. Dragovich, J. Fuchs, A. Gavrilik, H. Grosse, D.
Gurevich, J. Lukierski, M. Pavlov, H. Steinacker, Z. Raki\'c, W. Werner,
and S. Wo\-ro\-no\-wicz for fruitful discussions. Also he is grateful to
the Alexander von Humboldt Foundation for valuable support as well as to J.
Cuntz for kind hospitality at the Mathematisches Institut, Universit\"at
M\"unster, where this paper was finalized. Both authors would like to
express their gratitude to D. Shklyarov who attracted their attention to
the fact that the results of this work were actually valid for $q$ not
being a root of unit, rather than for $0<q<1$. We are also grateful to the
referee who pointed out some inconsistencies in a previous version of the
paper.

\newpage

\mbox{}

\bigskip

\textbf{Table 1.}

\bigskip

$$
\begin{array}{||c|c|c||}
\hline\hline
\textbf{Symbolic matrices} &
\mathbf{U_q(\mathfrak{sl}_2)}-\textbf{symmetries} &
\begin{array}{c}\textbf{Classical limit}\\
\mathbf{\mathfrak{sl}_2}-\textbf{actions} \\ \textbf{by differentiations}
\end{array}
\\ \hline\hline
\boldmath\left[
\begin{pmatrix}0 & 0\\ 0 & 0\end{pmatrix}_0;
\begin{pmatrix}0 & 0\\ 0 & 0\end{pmatrix}_1
\right] &
\begin{aligned}
k(x) &=\pm x,\;k(y)=\pm y, &
\\ e(x) &=e(y)=0, &
\\ f(x) &=f(y)=0, &
\end{aligned} &
\begin{aligned}
h(x) &=0,\quad h(y)=0, &
\\ e(x) &=e(y)=0, &
\\ f(x) &=f(y)=0, &
\end{aligned}
\\ \hline
\boldmath\left[
\begin{pmatrix}0 & \star\\ 0 & 0\end{pmatrix}_0;
\begin{pmatrix}0 & 0\\ 0 & 0\end{pmatrix}_1
\right] &
\begin{aligned}
k(x) &=qx, &
\\ k(y) &=q^{-2}y, &
\\ e(x) &=0,\quad e(y)=b_0, &
\\ f(x) &=b_0^{-1}xy, &
\\ f(y) &=-qb_0^{-1}y^2 &
\end{aligned} &
\begin{aligned}
h(x) &=x, &
\\ h(y) &=-2y, & \\
e(x) &=0,\quad e(y)=b_0, &
\\ f(x) &=b_0^{-1}xy, &
\\ f(y) &=-b_0^{-1}y^2 &
\end{aligned}
\\ \hline
\boldmath\left[
\begin{pmatrix}0 & 0\\ \star & 0\end{pmatrix}_0;
\begin{pmatrix}0 & 0\\ 0 & 0\end{pmatrix}_1
\right] &
\begin{aligned}
k(x) &=q^2x, &
\\ k(y) &=q^{-1}y, &
\\ e(x) &=-qc_0^{-1}x^2, &
\\ e(y) &=c_0^{-1}xy, &
\\ f(x) &=c_0,\quad f(y)=0, &
\end{aligned} &
\begin{aligned}
h(x) &=2x, &
\\ h(y) &=-y, &
\\ e(x) &=-c_0^{-1}x^2, &
\\ e(y) &=c_0^{-1}xy, &
\\ f(x) &=c_0,\quad f(y)=0. &
\end{aligned}
\\ \hline
\boldmath\left[
\begin{pmatrix}\star & 0\\ 0 & 0\end{pmatrix}_0;
\begin{pmatrix}0 & 0\\ 0 & 0\end{pmatrix}_1
\right] &
\begin{aligned}
k(x) &=q^{-2}x, &
\\ k(y) &=q^{-1}y, &
\\ e(x) &=a_0,\quad e(y)=0, &
\\ f(x) &=-qa_0^{-1}x^2+ty^4, &
\\ f(y) &=-qa_0^{-1}xy+sy^3. &
\end{aligned} &
\begin{aligned}
h(x) &=-2x, &
\\ h(y) &=-y, &
\\ e(x) &=a_0,\quad e(y)=0, &
\\ f(x) &=-a_0^{-1}x^2+ty^4, &
\\ f(y) &=-a_0^{-1}xy+sy^3. &
\end{aligned}
\\ \hline
\boldmath\left[
\begin{pmatrix}0 & 0\\ 0 & \star\end{pmatrix}_0;
\begin{pmatrix}0 & 0\\ 0 & 0\end{pmatrix}_1
\right] &
\begin{aligned}
k(x) &=qx,\;k(y)=q^2y, &
\\ e(x) &=-qd_0^{-1}xy+sx^3, &
\\ e(y) &=-qd_0^{-1}y^2+tx^4, &
\\ f(x) &=0,\quad f(y)=d_0, &
\end{aligned} &
\begin{aligned}
h(x) &=x,\quad h(y)=2y, &
\\ e(x) &=-d_0^{-1}xy+sx^3, &
\\ e(y) &=-d_0^{-1}y^2+tx^4, &
\\ f(x) &=0,\quad f(y)=d_0, &
\end{aligned}
\\ \hline
\boldmath\left[
\begin{pmatrix}0 & 0\\ 0 & 0\end{pmatrix}_0;
\begin{pmatrix}0 & \star\\ \star & 0\end{pmatrix}_1
\right] &
\begin{aligned}
k(x) &=qx, &
\\ k(y) &=q^{-1}y, &
\\ e(x) &=0,\quad e(y)=\tau x, &
\\ f(x) &=\tau^{-1}y,\;f(y)=0, &
\end{aligned} &
\begin{aligned}
h(x) &=x, &
\\ h(y) &=-y, &
\\ e(x) &=0,\quad e(y)=\tau x, &
\\ f(x) &=\tau^{-1}y,\;f(y)=0. &
\end{aligned}
\\ \hline\hline
\end{array}
$$

\newpage 

\mbox{}

\end{document}